\documentclass[12pt, a4paper]{article}
\usepackage[dvips]{graphicx}
\usepackage[mathscr]{eucal}
\usepackage{epsfig}
\usepackage{amssymb}
\renewcommand{\d}{\d}

\newcommand{\ov}{\overline}
\newcommand{\di}{\displaystyle}
\newcommand{\ctg}{\mathop{\rm cot}\nolimits}

\newcommand{\arccot}{\mathop{\rm arccot}\nolimits}

\usepackage[dvips]{color}
\begin{document}
\large
\begin{center}
{\bf\Large Prime number logarithmic geometry on the plane}\\[0.3cm]

{\large Lubomir  Alexandrov}\\[1cm]

{\it JINR, LTP, Dubna 141980, Moscow Region, Russia}\\[1cm]
\end{center}

\begin{quote}
{\bf Abstract}
{\it
We found a regularity of the behavior of primes that allows to represent both
prime and natural numbers as infinite matrices with a common formation rule
of their rows. This regularity determines a new class of infinite cyclic groups
that permit the proposition a plane--spiral geometric
concept of the arithmetic.}\\

\end{quote}

\section{Introduction}

Counting arithmetic functions for different prime sets can be assigned
to {\it the archaic mathematical reality}.

Nevertheless, the generated by them prime sequences, named {\it
Eratosthenes progressions},  became known only in recent years
(e.g.,\cite{b1},\cite{b2},\cite{b3} and sequences A007097, A063502, A064110 in \cite{b4}).

 The Eratosthenes progression possesses a common formation law
of its elements (an {\it inner prime number distribution law})
the realization of which is based on a multiple use of the Eratosthenes
sieve \cite{b1} (Figure 1).

The derivation of Eratosthenes progressions and their systematic
investigation is directed to a learning the nonasymptotic  behaviour of primes,
i.e., of the function's behaviour
$$
d (n)=p(n+1)-p(n),\, n=1,2,\ldots, \ov{n},
$$
where $p(n)$ is the $n$th prime and
$\ov{n}$ is a sufficiently large natural number.

 The inner prime number distribution  law
can be applied mostly in mathematics
itself, for example, when constructing new geometric concepts in
arithmetic.

Following Alain Connes (\cite{b5} pp. 208--209), it can be supposed
that the specific behaviour of primes will reflect
itself in the new geometry sought for understanding quantum gravity.

In biochemistry, the specific behaviour of primes can manifest itself
in the laws of formation and functioning of large molecules, from
$10^3$--atomic insulin and hemoglobin up to $3\cdot 10^5$--atomic
proteins and enzymes.

In this paper, the general statement of the problem for
derivation of Eratosthenes progressions is given and their
basic properties are presented. The general results are applied
to the sequence of primes itself
$$
P=\{ 2, 3, 5, 7, 11,\ldots\} =\{p(n)\}_{n=1, 2,\ldots},
$$
as well as to the following related to $P$ sequences:\\

\noindent
$M=\mathbb{N}\setminus P=\{4, 6, 8, 9,\ldots\}=
\{m(n)\}_{n=1, 2,\ldots}$ the set of composite numbers;\\

\noindent
$T=\{t(\nu)=(p(\nu),p(\nu+1)):\nu\in \Lambda\}$
the set of twin pairs, where \linebreak\phantom{ffffffff}
$\Lambda =\{n:p(n+1)-p(n)=2, n\in \mathbb{N}\}=
\{2,3,5,7,10,13,17,\ldots\}$;\\

\noindent
$T_1=\{p(\nu):(p(\nu),p(\nu+1))\in T,\,\nu\in \Lambda\}=
\{t_1(\nu)\}_{\nu\in\Lambda}$ the set\linebreak\phantom{ffffffffff}
 of first elements of twins;\\

\noindent
$T_2=\{p(\nu+1):(p(\nu),p(\nu+1))\in T,\,\nu\in\Lambda\}=
\{t_2(n)\}_{\nu\in\Lambda}$ the\linebreak\phantom{ffffffffff}
 set of second elements of twins;\\

\noindent
$T_{3}=T_1\cup T_2=\{3,5,7,11,13,17,19,\ldots\}$  the set of
twin elements;\\

\noindent
$S=P\setminus T_{3}=\{2,23,37,47,53,\ldots\}$  the set of
isolated primes\cite{b4}, A007510;\\

\noindent
$D_{6n-1}=
\{6n-1\in P: n=1, 2,\ldots,\}=\{5, 11, 17,\ldots\}$  the
 set of primes of\linebreak\phantom{ddfffffffff} the kind $6n-1$;\\

\noindent
$D_{6n+1}=\{6n+1\in P: n=1, 2,\ldots,\}=\{7, 13, 19,\ldots\}$  the
 set of primes of\linebreak\phantom{ddfffffffff} the kind  $6n+1$, and \\

\begin{figure}
\centering
\includegraphics[width=16cm, height=20cm]{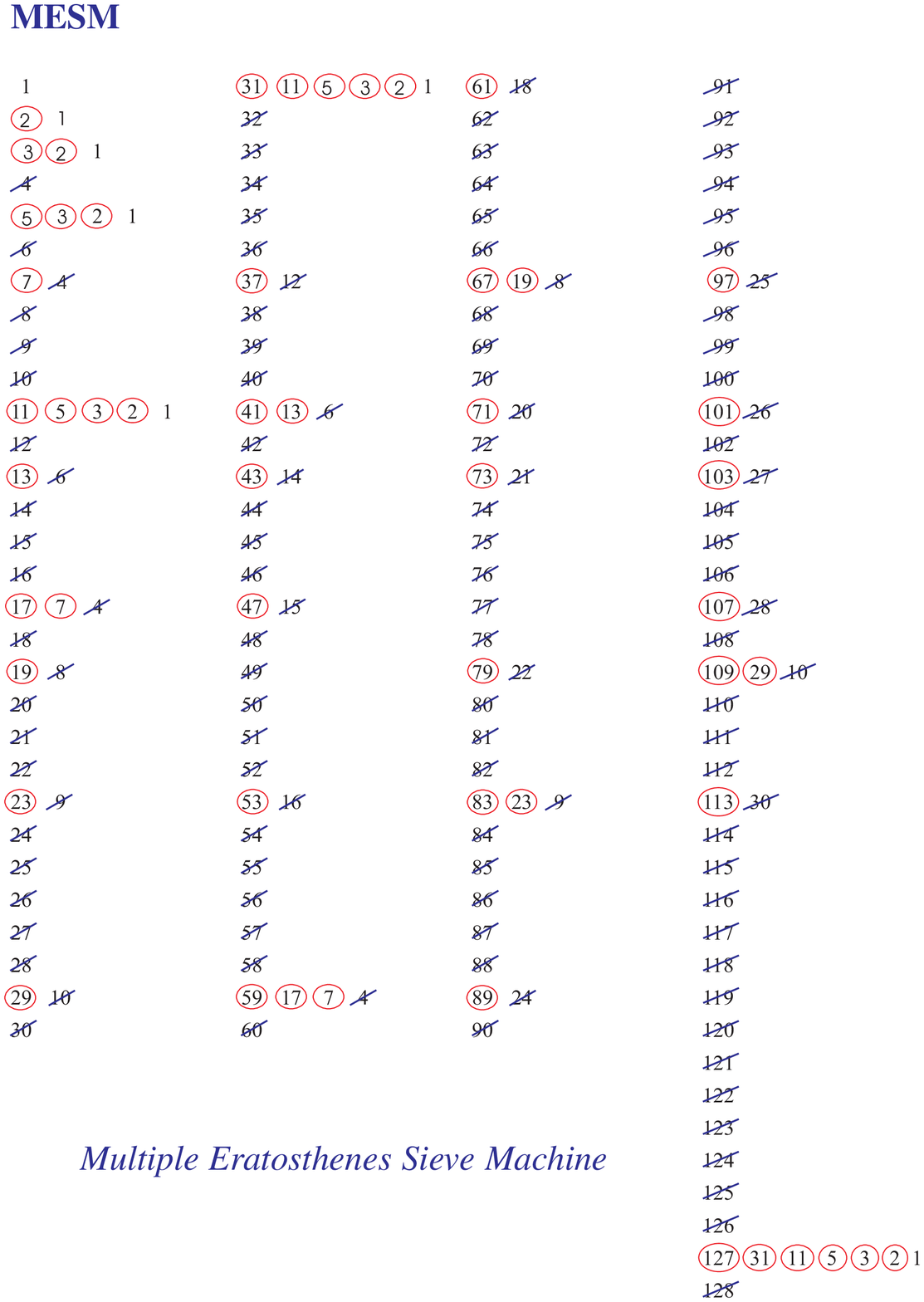}
\caption{}
\label{fig1}
\end{figure}

\noindent
$T_4=\{t(n): (t_1(n)+t_2(n))/2 =6\cdot q,\, q\in P\}$  the set of
twins with\linebreak\phantom{fffffffff} minimal average
(\cite{b2}, p. 15).\\

The sets $T,\; T_1 - T_4$\mbox{ and } S below will be supposed to be infinite.

In this paper some properties of Eratosthenes progression such as
distribution laws of the progression elements,
$\zeta$--functions for the progressions and their connection
with the Riemann  $\zeta$--function are only mentioned.

The main result of this paper consists in the proposed plane--spiral
geometric concept of arithmetic, compatible with
the linear Cartesian concept.

The real semiaxis  $\mathbb{R}_{+}^{1}$ in the new geometric model
is isometrically mapped as {\it a logarithmic  spline-spiral} on the
plane $\mathbb{R}^{2}$ in such a way that the Eratosthenes rays,
not intersecting each other, cross the spiral only at the primes.

The spiral arithmetic allows one to interpret in a new way
the basic counting function  $\pi (x)$,\,
the Littlwood's $\Omega$--theorem
and also  gives an arithmetic interpretation of the
distribution in natural series of all kinds of  clusters of primes
(see \cite{b6}, for example) and twin pairs, in particular.

The basic object in the spiral geometry is a {\it spider--web}
$W_n$ composed of spiral and Eratosthenes rays intersecting it,
in which the number of rotations $n$ infinitely increases.

 The web $W_n$ consists of embedded
{\it concave--convex trapezoids} of primes with a
characteristic formation law. This law is a direct consequence of the
inner prime number distribution law.

The plane $\mathbb{R}^{2}$ is considered as  {\it a mosaic}
composed of {\it elementary\linebreak concave--convex trapezoids}.

The web $W_n$ {\it geometrically select (personalyzes)} primes, and also all kinds
of linear and plane configurations  of primes.

\section{Splitting theorem for infinite sequences of primes}
\subsection{Basic definitions}

Let sets  $A\subset\mathbb{N}$ and
$B\subset\mathbb{N}$ with the properties
\begin{equation}
\label{eq1}
A\cap B=\varnothing,
\end{equation}
\begin{equation}
\label{eq2}
A\cup \ov{B}=\mathbb{N},
\end{equation}
where $\ov{B}=\{1\}\cup B$ are given.

Let the arithmetic function
$$
g(n): \mathbb{N}\to A
$$
generate (denote) the $n$th element $a(n)\in A$.

Then  the {\it counting recurrent law }
\begin{equation}
\label{eq3}
\varepsilon_{a(0)}^+: a(n+1)=g(a(n)),\, n=0, 1, 2,\ldots, a(0)\in\mathbb{N}
\end{equation}
determines {\it an A--counting progression  $\varepsilon_{a(0)}^+$} and
{\it an $A$--counting ray}
$$
r_{a(0)}=\{a(n): a(n+1)=g(a(n)),\, n=0, 1, 2,\ldots, a(0)\in
\mathbb{N} \}.
$$

Together with the function  $g(n)$, its inverse function, the {\it
$n$th number of element} $a(n)\in A$, is also uniquely
determined (in a purely arithmetical sense it is a counting function)
$$
g_{-1}(a): A\to \mathbb{N}.
$$

The functions $g(n)$ and $g_{-1}(a)$ are strictly monotonic and satisfy
the equalities
$$
g(g_{-1}(a))=a,\quad  g_{-1}(g(n))=n.
$$

By means of $g(n)$ and $g_{-1}(a)$ the compositions
$$
\begin{array}{l}
g_{n}(a(0))=\underbrace{g(\ldots g}_{n}(a(0))\ldots ) =a(n), \\
\phantom{rrr}\\
g_{-n}(a(n_1))=
\underbrace{g_{-1}(\ldots g_{-1}}_{n} (a(n_1))\ldots ),\,
\quad\mbox{with }n\leq n_1.
\end{array}
$$
are introduced.

These compositions satisfy  the equalities
$$
\begin{array}{l}
g_{n_1}(g_{n_2}(a(0)))=g_{n_1+n_2}(a(0)),\quad n_1,
n_2\geq 1, \\
\phantom{ff}\\
g_{-n_1}(g_{n_2}(a(0)))=g_{n_2-n_1}(a(0)),\quad
1\leq n_1\leq n_2.
\end{array}
$$

An extension of the $A$--counting progression
 $\varepsilon_{a(0)}^+$
with negative numbers $\varepsilon_{a(0)}^{-}=-\varepsilon_{a(0)}^+$
leads to an infinite cyclic group
\begin{equation}
\label{group}
\varepsilon_{a(0)}=\varepsilon_{a(0)}^{-}\cup\{a(0)\}\cup
\varepsilon_{a(0)}^+,\quad g_{-n}(a(0))=-g_n(a(0)),\; n>0
\end{equation}
under composition $g_n(a(0)),$ with {\it a depth} $n\in\mathbb{Z}$ and
a {\it generator} $a(0)\in\ov{B}$.

Two elements from  $\varepsilon_{a(0)}$ interact under the composition
rule
\begin{equation}
\label{group1}
g_{n_1}(a(0))\circ g_{n_2}(a(0))=g_{n_1}(g_{n_2}(a(0)))=
g_{n_1+n_2}(a(0)),\;
 n_1, n_2 \in\mathbb{Z}.
\end{equation}

\subsection{Basic assertion and its consequences}

The following assertion about the splitting of the set
$A$ in a denumerable number
of denumerable subsets with a common law(\ref{eq3}) of formation
of its elements is given:\\

{\bf Theorem 1.}\,
{\it For any sets $A$ and $B$ with properties  (\ref{eq1}) and
(\ref{eq2}) the following equalities hold
$$
\begin{array}{l}
\bigcap\limits_{a(0)\in \ov{B}}r_{a(0)}=\varnothing,
\phantom{eee}
\bigcup\limits_{a(0)\in \ov{B}}r_{a(0)}=A\;.\\
\end{array}
$$}

Theorem 1 leads to a matrix representation of the sequences
$A$ and $\mathbb{N}$ with peculiar properties of their elements.\\

{\bf Corollary 1.\,}{\it There exists an one-to-one mapping}
\begin{equation}
\label{eq4}
\ov{\varphi} (a(0)): \ov{B}\to\phantom{,} ^2A=\{r_{a(0)}\}_{a(0)\in\ov{B}}\,
\equiv \{a_{\mu\nu}\}_{\mu , \nu =1, 2,\ldots}
\end{equation}
( $^2A$ denotes the matrix representation of the elements of $A$).\\

From (\ref{eq4})  a matrix representation to the natural series
\begin{equation}
\label{eq5}
^2\mathbb{N}=\|\ov{B}\phantom{d,} ^2A\|,
\end{equation}
 where $\ov{B}=\mbox{Column}\{a_{\mu 0}\}_{\mu =1, 2,\ldots}$ also
 follows.

The matrices  $^2 A$ and $^2\mathbb{N}$ shall be called
{\it mesm}--matrices.

In the case when  $A=P$ and $B=M$ an example of the left upper corner
of the matrix  $^2\mathbb{N}$ (\cite{b2}, pp. 18--22)
is given in Appendix 1.\\[0.5cm]

{\bf Corollary 2.}\,
{\it The rows of matrices  $^2\mathbb{N}$ are isomorphic to the
row $r_{1}$ with respect to the mapping}
$$
\Psi(g_{n}(1)): g_{n}(1)\to g_{-n}(g_{n}(1))\to a(0)\to g_{n}(a(0)),
\quad a(0)\in \ov{B}, \quad a(0)>1.
$$
{\it The columns of the matrix  $^2A$ are isomorphic to the column $\ov{B}$
with respect to the mapping}
$$
\varphi(a(0)): a(0)\to g_{n}(a(0)),\quad a(0)\in \ov{B}.
$$

\vspace{0.5cm}

In the case  $A=P$ and $B=M$ Figure 2 illustrates mentioned isomorphisms.
In Figure 2, an one-to-one correspondence between {\it rooted trees} and
elements of $\mathbb{N}$, proposed by F. G$\ddot{\mbox{o}}$bel \cite{b7}
is used (see the $1th$ row of Figure 2).

Theorem 1 leads also to an important consequence, which reveals the
arithmetic
nature of the fine structure of the set $A$ elements' distribution  among
the natural numbers.

Let $g_{-1}(n', n''),\; n', n'' \in\mathbb{N}$ denote the number
of elements $A$ in the interval  $(n', n'')$.

{\bf Corollary 3.}\,
{\it For the matrix $[\ov{B}\phantom{d} ^2A]$ elements
the following equalities hold:}
\begin{equation}
\label{eq6}
\left.
\begin{array}{l}
g_{-1}(a_{\mu 0}, 0)=a_{\mu_1}-1, \quad \mu =1, 2,\ldots,\\
\phantom{ffff}\\
g_{-1}(a_{\mu_1 \nu_1}, a_{\mu_2 \nu_2})=|a_{\mu_1 (\nu_1-1)}-
a_{\mu_2 (\nu_2-1)}|-1,\\
\phantom{ddddd}\\
\mu_{i}, \nu_{i}\geq1, \quad i=1, 2 .
\end{array}
\right\}
\end{equation}
\section{The Theorem 1 application to special
cases of sets $\mathbf{A}$ and $\mathbf{B}$}
\subsection{About new $\mathbf{A}$--counting progressions}

In the case when  $A$ and $B$ take usual values the law
(\ref{eq3}) generates known  $A$--counting progressions.
So, for example, at $A=\{even \}$ and $B=\{odd \}$ a generating function
is of kind  $g(a(n))=2a(n)-1$ and in this case
$\varepsilon_{2}^+=\{2, 3, 5, 9, 17, 33, 65, 129\ldots\}$
is a {\it Pisot sequence} (\cite{b4}, A000051).

New $A$--progressions one occur when the behaviour of $A$ elements
among natural numbers is unknown and it cannot be considered as a
probabilistic. Besides the sequences of primes $P$, all subsequences of
$P$,  in the formation of which the Eratosthenes sieve combines with
an additional deterministic filter $f(n)$ (this is the formation rule
of the considered subsequence), should also be considered belonging
to this class. The set of these subsequences
shall be denoted by  $\mathscr{E}_f$.

The Dirichlet theorem about the existence of infinite primes of the kind
$\alpha n +\beta$ (an additional filter) for arbitrary coprimes
$\alpha$ and $\beta$ shows that  $\mathscr{E}_f$ is infinite.

In particular, we have inclusions
$T_1, T_2, T_{3}, S\in\mathscr{E}_f \quad and \quad
D_{\alpha n\pm 1}\in\mathscr{E}_f\ at\;\alpha =4,6$.

For all elements $\mathscr{E}_f$ there exists a {\it mesm$_f$}--process,
which is analogous to the process represented in Figure 1.
From the  $A$--split theorem
it follows that for every $A_f\in\mathscr{E}_f$ and
$B_f=\mathbb{N}\setminus A_f$ there exists a {\it mesm$_f$}--matrix
$[B_f\; ^2A_f].$

As a result of a {\it mesm}--transition $A_f\to\phantom{,} ^2A_f$,
the elements of the rows
$^2P,\, ^2T_1,\, ^2T_2,\, ^2T_{3},\, ^2S$
\mbox{ and} $^2D_{\alpha n\pm 1 }(\alpha =4,6)$ already will be
distributed according to the {\it inner law} (\ref{eq3}), which now should
be understood as a specific {\it self-smoothing } (only with respect to the
rows $^2A_f$) of the irregularities in the appearance of the
elements $A_f$ in the natural series.

\subsection{The basic case: $\mathbf{A=P}$ and $\mathbf{B=M}$.}

The upper left corner of the matrix  $^2P$ and its extension to the
matrix $^2\mathbb{N}$ are represented in Appendix 1.
The Theorem 1 has been proved inductively in \cite{b2}, pp. 4--8.

The first elements of the first rows of the matrix $^2P$
were primarily determined {\it by hand}
by means of {\it MESM} (Figure 1).
In such a way the law (\ref{eq3}) with $g(n)=p(n)$ ({\it
Eratosthenes progressions})  was discovered
\cite{b1}.

The extension of the matrix $^2P$ rows on negative primes according
to the rules (\ref{group}), (\ref{group1})
leads to infinite cyclic groups under composition
$p_n(a(0)),$\linebreak $n\in\mathbb{Z}$
with a generators $a(0)\in M$. An example of such a group
is the set
$$
\varepsilon_4 =\{\ldots,-p_n(4),\ldots,-59, -17,-7,4,7,17,59,\ldots,
p_n(4),\ldots\}.
$$

A part of $^2P$ represented in Appendix 1 has been
computed by means of {\bf Mathematica} function
{\bf \mbox{NestList[Prime, a(0), n]}}.

\begin{figure}
\centering
\includegraphics[width=16cm, height=20cm]{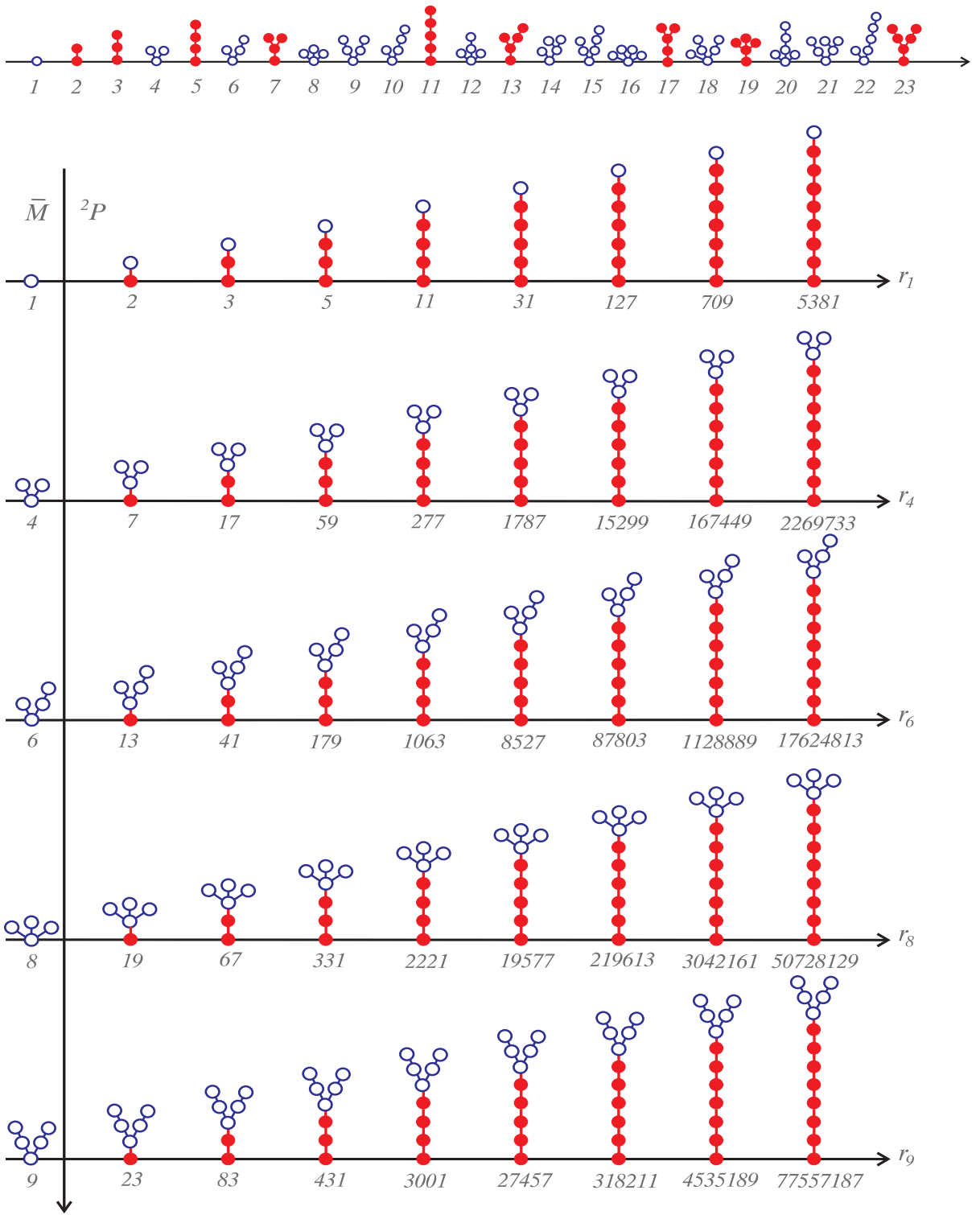}
\caption{{\it "MESM {\rm\&} F. G$\ddot{\mbox{o}}$bel" forest of
rooted trees }}
\label{fig2}
\end{figure}

The row elements of  the matrix $[\ov{M}\phantom{d}\, ^2P]$ determine
new subsets of natural numbers
$$
N_m=\{p_{n_1}^{\alpha_1}(m)\ldots p_{n_k}^{\alpha_k}(m):\, \forall\, n_i,
\alpha_i \in\mathbb{N},\, i=1,2,\ldots,k,\, \forall\, k\in\mathbb{N}\},\,
m\in \ov{M}:
$$
$$
N_1=\{2,3, 2^2, 5,2\cdot 3, 2^3, 3^2, 2\cdot 5,11, 2^2\cdot 3, 3\cdot 5,
2^4,\ldots \},
$$
$$
N_4=\{7, 17, 7^2, 59, 7\cdot 17, 277, 17^2, 7^3,
7\cdot 59, 7^2\cdot 17, 1787, 7^4,
\ldots \},
$$
$$
N_6=\{13, 41, 13^2, 179, 13\cdot 41, 1063, 41^2, 13^3, 13^2\cdot 41, \ldots
\},
$$
and so on.

According to Corollary 2, the behaviour of composite numbers
reflects on the behaviour of the elements of the columns of the matrix
$^2P$.

On the other hand, the structure of the set $M$
depends on  the structure of the set
of primes because $M$ can be represented as a chain of
$\alpha_{\mu}$--element segments
$(\mbox{ where }\alpha_{\mu}=d (\mu )-1)$
from consequent composite numbers
$$
\ov{m}_{\mu}(\alpha_{\mu})=
\{p(\mu)+1, p(\mu)+2,\ldots, p(\mu +1)-1\},\quad \mu =2, 3,\ldots
$$
$$
( \ov{m}_{2}(1)=\{4\},\, \ov{m}_{3}(1)=\{6\},\, \ov{m}_{4}(3)=
\{8, 9, 10\},\ldots ).
$$

The segments are connected in a whole set $M$ by means of {\it ghost
primes} $\omega_{\mu}=\langle p(\mu)\rangle$\quad
$(\omega_{2}=\langle 3\rangle,\, \omega_{3}=\langle 5\rangle,\, \omega_{4}=
\langle 7\rangle,\ldots )$.

The Eratosthenes  progressions
$\{\varepsilon_{m}^+\}_{m \in\ov{M}}$ (i.e., rows of the matrix
$^2P$) conform to the inner prime number distribution law
\begin{equation}
\label{eqnew1}
a(n+1)=p(a(n))=p_{n+1} (a(0)),\quad n=0, 1, 2,\ldots , a(0)\equiv m\in\ov{M},
\end{equation}
but the deviation of the rows $^2P$ between each other (i.e., the
distribution of primes in the columns of $^2P$) again persists
dependent of the oddish behaviour of primes.

The main information left out of the inner law (\ref{eqnew1})  is
reflected in the structure of the first matrix $^2P$ column
$$
P_1=column[p_{11}, p_{21},\ldots, p_{\mu 1},\ldots].
$$
The following assertion about the $P_1$ structure
is valid.\\

{\bf Theorem 2.}\,{\it Mapping
$\varphi : a(0)\to p_{\mu 1}$ defines a correspondence between segments
of composite numbers  $\ov{m}_{\mu} (\alpha_{\mu})$ and clusters of
 $\alpha_{\mu}$--successive primes
$$
c_{\mu}(\alpha_{\mu})=\{ p_1(p(\mu)+1), p_1(p(\mu)+2),\ldots, p_1
(p(\mu +1)-1) \}\subset P
$$
 in the cases $\alpha_{\mu}\geq 3$, and separate primes
$p_1(p(\mu)+1)$ in the cases $\alpha_{\mu}=1$. At their ends the clusters
are complemented by the ghost images up to prime number segments
$$
\ov{c}_{\mu}(\alpha_{\mu})=\{ p_1(\langle p(\mu)\rangle),
c_{\mu}(\alpha_{\mu}), p_1(\langle
p(\mu+1)\rangle)\}
$$
and the equality $P=\bigcup\limits_{\mu=1}^{\infty}
\ov{c}_{\mu}(\alpha_{\mu})$ is fulfilled.}

The next theorem about twin pairs
$t({\nu})=\left(t_1 ({\nu}),\, t_2({\nu})\right)\in T,
\quad {\nu}=3,5,7,\ldots$
 is also justified.\\[0.5cm]

{\bf Theorem 3.\,}
{\it For each pair $t({\nu})$(after the pair (3,5))
at least one of the elements $t_1({\nu})$
or $t_2({\nu})$ belongs to the first  column $P_1$ .

The mapping $\varphi^{-1}:\, p_{\mu 1}\to m_{\mu}$ defines a correspondence
between pairs with both elements on $P_1$ ($u$--twin) and pairs of
subsequent elements of some segment
$\ov{m}_{\ov{\mu}}(\alpha_{\ov{\mu}})\subset M$ with
$\alpha_{\ov{\mu}}\geq 3$.

For a pair  with one element $t_1({\nu})$ $(or\; t_2({\nu}))$ on $P_1$
($b$--twin) the mapping\linebreak $\varphi^{-1}:\, p_{\mu 1}\to m (\mu)$
associates $t_1({\nu})$ ($or\; t_2({\nu}))$ with the element  $p(n)+1$, or
the element  $p(n+1)-1$ of some segment $\ov{m}_{\ov{\mu}}
(\alpha_{\ov{\mu}})\subset M$ at $\alpha_{\ov{\mu}}\geq 3$, or with the
element of some one-element segment $\ov{m}_{\ov{\mu}}(1)\subset M$.

The mapping $\varphi^{-1}:\, p_{\mu_1 \nu_1}\to$ $p_{\mu_1 (\nu_1 -1)},\,
\nu_1\geq 2$ relates the second element
 $t_2({\nu})$ ($or\; t_1({\nu}))$ to one of the ghosts
$\langle p(\mu )\rangle\equiv p_{\mu_1 (\nu_1 -1)}$ or
$\langle p(\mu +1)\rangle\equiv$ $p_{\mu_1 (\nu_1 -1)}$.}\\[0.5cm]

The following properties of the matrix $^2P$ rows and columns
are briefly veiwed:

\begin{itemize}
\item[$\mathbf{q_1)}$] The difference
$d_{m}(n)=p_{(n+1)}(m)-p_{n}(m),\quad  n=1,2,
\ldots,\, m\in\ov{M}$
monotonically increase under the estimate
$$
d_{m} (n)>p_{n}(m)(\ln p_{n}(m)-1)
$$
unlike the difference  $d(n)$ whose  behaviour
only on the face of it may seems to be a chaotical one \cite{b8};

\item[$\mathbf{q_2)}$] The sequence
$\eta (s,\, m)=\sum\limits_{n=1}^{\infty}\frac{\di 1}{\di p^s_n (m)}$
converge for all
$m\in\ov{M}$ and\linebreak $s\geq 1.$

Note especially the convergence of the sum $\eta (1,\, m)$
(\cite{b2}, p. 10) when the sum
$\sum\limits_{n=1}^{\infty}\frac{\di{1}}{\di{p(n)}}$  diverges;

\item[$\mathbf{q_3)}$]  An analogue of the Euler  identity exists

$$
\zeta (s, m))\equiv 1+\sum\limits_{n\in N_{m}}\frac{1}{n^s}=
\prod\limits_{n=1}^{\infty}\left( 1- \frac{1}{p_{n}^{s}(m)}
\right)^{-1},\; m\in\ov{M}, s\geq 1;
$$

\item[$\mathbf{q_4)}$] The Riemann function
$\zeta (s)=
\sum\limits_{n=1}^{\infty}\frac{\di{ 1}}{\di{n^{s}}},\quad s
\in\, \not\hspace{-0.15cm}C $
can be represented by the functions $\zeta (s,\, m)$
$$
\zeta (s)=\prod\limits_{m\in\ov{M}}\zeta(s,\, m);
$$

\item[$\mathbf{q_5)}$] The asymptotic law for the primes and the simplified
Riemann formula for $\pi (x)$ give an opportunity to find approximately
$p_{n+1}(m)$,\linebreak $m\in\ov{M}$ by solving the equations with respect to $x$
\begin{equation}
\label{eq88}
L(x)=p_{n}(m),
\end{equation}
\begin{equation}
\label{eq99}
R(x)=p_{n}(m),
\end{equation}
$$
\mbox{where }\, L(x)=\int\limits_{0}^{x}\frac{\di ds}{\di\ln (s)},\;
R(x)=\sum\limits_{k=1}^{\infty}\frac{\di\mu (k)}{\di k} L (x^{1/k})
$$
and $\mu (k)$  is a M$\ddot{\mbox{o}}$bius function;

\item[$\mathbf{q_6)}$]  There exists an  approximate formula
\begin{equation}
\label{eq7}
\di n=\int\limits_{\alpha}^{p_n(\beta)}\frac{ds}{s\ln\ln s}+\varepsilon (n,
\,\beta),
\end{equation}
$\mbox{where } \alpha =11,\, \beta =1,\;
  n>4 \mbox{ for } r_1,\;
 \alpha =7,\, \beta =4 \mbox{ for } r_4\;
 \mbox{and }$
 $\alpha = \beta =m \mbox{ for the other }
 \mbox{rays } r_m.$

The absolute error $|\varepsilon|$ for the part of the matrix
$[\ov{M}\phantom{d} ^2P]$ in Appendix 1 is not greater than
$0.2$ when $n$ is small and $0.06$ when $n$ is large.

 Formula (\ref{eq7}) is a {\it prime number distribution law}
 of the rays $^2P$.

On Figure 3, the behaviour of the function (\ref{eq7}) is presented for
the ray $r_9$;

\item[$\mathbf{q_7)}$]  It is obvious that for the number
$\mu$ of the element $p_{\mu n}$ in the matrix $^2P$ column
$$
P_n =colomn [p_{1n}, p_{2n},\ldots,  p_{\mu n},\ldots]
$$
there exists an asymptotic formula
\begin{equation}
\label{eq8}
\mu\sim m-\int\limits_{2}^{m}\frac{ds}{\ln s}.
\end{equation}

This is the column $^2P$ {\it prime number distribution law}.

In order to use (\ref{eq7}) and (\ref{eq8}) it is necessary to know
the composite number $m$.
\end{itemize}

\subsection{ About other $\mathbf{A}$--counting progressions}

{\normalsize Applying the $A$--split theorem in the cases
$$
\begin{array}{l}
\phantom{eee}A=T_1\mbox{ and }B_1=T_2\cup M,\\
\phantom{eee}\\
\phantom{eee}A=S\mbox{ and }B_2=M\cup T_3,\\
\phantom{eee}\\
\phantom{eee}A=D_{6n-1}\mbox{ and }B_3=M\cup D_{6n+1}\cup\{2, 3\},\\
\mbox{and}\\
\phantom{eee}A=D_{6n+1},\, B_4=M\cup D_{6n-1}\cup\{2, 3\},\\
\end{array}
$$
we can obtain the next {\it mesm}-matrices of kind (\ref{eq5}):
$$
[B_1\; ^2T_1]=\left[
\begin{array}{lllllll}
1& 3& 11& 137& 5639& 641129& 152921807\ldots \\
2& 5& 29& 641& 44381& 7212059& \ldots\\
4& 17& 239& 12161& 1583927& \ldots& \\
6& 41& 1151& 93251&  16989317& \ldots& \\
7& 59& 1931& 176021& 35263691& \ldots& \\
8& 71& 2339& 221201& 45749309& \ldots& \\
\centerdot& \centerdot& \centerdot& \centerdot& \centerdot& \ldots&
\end{array}
\right];
$$

\vspace{1cm}

\begin{figure}
\centering
\includegraphics[width=12.5cm, height=8.5cm]{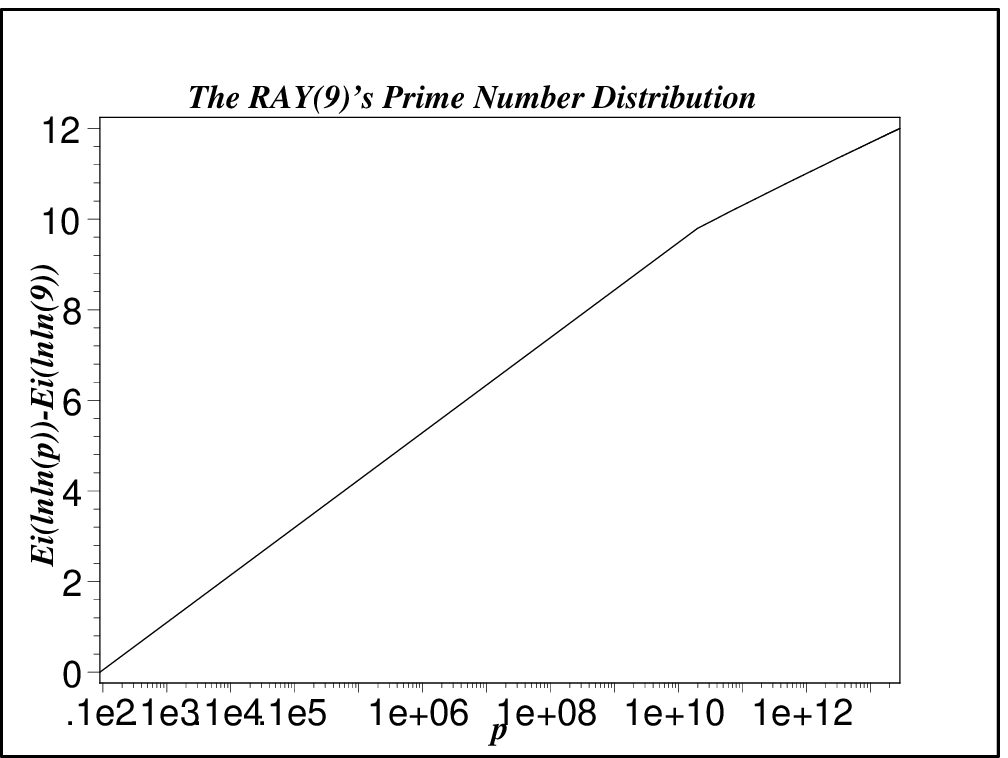}
\caption{}
\label{fig2}
\end{figure}
$$
[B_2\; ^2S]=\left[
\begin{array}{llllllll}
1& 2& 23& 263& 2917& 38639& 603311& 11093633\ldots \\
3& 37& 397& 4751& 64403& 1038629& 19661749& \ldots\\
4& 47& 491& 5897& 81131& 1328167& 25467419& \ldots \\
5& 53& 557& 6709&  93287& 1541191& 29778547& \ldots\\
\centerdot& \centerdot& \centerdot& \centerdot& \centerdot& \centerdot&
\centerdot& \ldots\\
22& 257& 2861& 37799& 589181& 10821757& 230452837&
\ldots \\
24& 277& 3079& 40823& 640121& 11807167& 252480587&
\ldots \\
\centerdot& \centerdot& \centerdot& \centerdot& \centerdot& \centerdot&
\centerdot& \ldots\\
\end{array}
\right];
$$

\vspace{1cm}
$$
[B_3\; ^2D_{6n-1}]=\left[
\begin{array}{llllllll}
1& 5& 29& 263& 3767& 76253& 2049263& 69633521\ldots \\
2& 11& 83& 953& 16223& 381221& 11579489& \ldots\\
3& 17& 137& 1721& 31883& 795803& 25434641& \ldots \\
4& 23& 197& 2663&  51803& 1348961& 44635001& \ldots\\
6& 41& 419& 6329& 135347& 3808109& 134441441& \ldots \\
\centerdot& \centerdot& \centerdot& \centerdot& \centerdot& \centerdot&
\centerdot& \ldots\\
\end{array}
\right];
$$

\vspace{1cm}
$$
[B_4\; ^2D_{6n+1}]=\left[
\begin{array}{llllllll}
1& 7& 61& 727& 12343& 284083& 8457367& 312953941\ldots \\
2& 13& 109& 1429& 26113& 642937& 20262883&
787318099\ldots\\
3& 19& 181& 2539& 49669& 1291471& 42627997& \ldots \\
4& 31& 331& 5011&  105277& 2908753& 10144807& \ldots\\
5& 37& 397& 6211& 133633& 3761239& 132710947&
\ldots \\
\centerdot& \centerdot& \centerdot& \centerdot& \centerdot& \centerdot&
\centerdot& \ldots\\
\end{array}
\right].
$$

The matrices $[B_1\; ^2T_{1}]$, $[B_2\; ^2S]$ and $[P\phantom{d} ^2M]$
were published in \cite{b4} as $A063502$, $A064110$ and $A025003$--$A025006$,
respectively. Matrices $[B_3\; ^2D_{6n-1}]$ and $[B_4\; ^2D_{6n+1}]$ are
the new ones.

New {\it mesm}--matrices can be obtained also for the Euler primes
of the kind\linebreak $n^2+n+41$ ($r_1=\{41, 1847, 1573316,\ldots\}$), and for the
Hardy-Littlwood primes of the kind $H_{n^2 +1}=\{n^2+1\in P:
n=1,2,\ldots\}$ where at $B_5=\mathbb{N}\setminus H_{n^2+1}$ we have

$$
[B_5\; ^2H_{n^2+1}]=\left[
\begin{array}{llllll}
1& 2& 5& 101& 746497&  \phantom{v}286961228404901\ldots \\
3& 17& 7057& 11424189457\ldots\\
4& 37& 44101& 637723627777\ldots \\
6& 197& 3496901\ldots\\
7& 257& 6421157\ldots \\
\centerdot& \centerdot& \centerdot& \centerdot\ldots\\
\end{array}
\right].
$$}

All pointed out  {\it mesm$_f$}--matrices are not studied.
In particular, an analogue  of the  distribution laws  (\ref{eq7}) and
(\ref{eq8}) has not been found for them with the exception of the matrix
$[P\phantom{d} ^2M]$ for which an analogue of the law  (\ref{eq8})
is known. However,  the common Corollaries 2 and 3
of the $A$--split theorem remain valid for them.

\section{Logarithmic geometry of primes on the plane.}
\subsection{The Prime Number Spider Web  (PNSW) Hypothesis}

One of the main application of {\it the prime number distribution law}
 (\ref{eqnew1})
consists in constructing the plane spiral geometric concept  of arithmetic.

Let
$$
\mathscr{L}_f=\{\rho (\theta)= {\left(f(\theta)\right)}^{\di\theta}:\,
f(\theta)\in C^1[0,\infty),\,
f(\theta)\geq 1,\, 0\leq \theta <\infty \}
$$
denote a class of logarithmic spirals with an arc length
$$
\lambda (0, \theta)=\int\limits_{\di 0}^{\di\theta} \left(f(x)\right)^{\di x}\left(
\left( \ln f(x)+\frac{xf'(x)}{f(x)}\right)^2+1\right)^{1/2}dx.
$$
The plane spiral geometric concept is based
on the following PNSW--hypothesis
\cite{b1}.\\

{\bf Conjecture  1.\,}{\it
On the plane $\mathbb{R}^2$ there exists a unique spiral
$\ov{\rho}(\theta)\in\mathscr{L}_f$ and the  corresponding to it sets of
angles
$$
\{\theta_{mn}\}_{m\in\ov{M}}\, ,\, n=1,2,\ldots,\quad\theta_{mn'}
<\theta_{mn''}\mbox{ at } n'<n''
$$
such that the following conditions are fulfilled:\\

(i)\quad $\lambda (0, \theta_{mn})=p_n(m),\, n=1,2,\ldots,\;
m\in\ov{M} ;$\\

(ii) the primes $p_n(m),\, n=1,2,\ldots$ lie on  the same ray
$\ell_{m}\subset R^2$ with a positive direction corresponding to
increasing $n$;\\

(iii) two arbitrary rays $\ell_{m_1}$ and
$\ell_{m_2},\; m_1, m_2\in\ov{M}$ do not intersect
each other and are non--parallel.}

\subsection{Logarithmic spline--spiral}

Under the substitution $f(\theta)=e^{\di \ctg\varphi},\,\mathscr{L}_f$
turns in a one--parametric family of logarithmic spirals
$$
\mathscr{L}_{\varphi}=\left\{\rho_{\varphi}=e^{\di (\ctg\varphi )\theta}:\,
0<\varphi<\frac{\pi}{2},\, 0\leq \theta<\infty\right\}
$$
with an arc length
\begin{equation}
\label{eqnew2}
\lambda (0, \theta)=\frac{1}{\cos\varphi}\left( e^{\di (\ctg\varphi)\theta}-1\right).
\end{equation}
Now the required by the  PNSW--hypothesis sets of angles with
respect to $m$ and $n$, according to the condition (i),
are given by the formula
$$
\theta_{mn}=\tan \varphi\ln (p_n(m)\cos\varphi +1)).
$$

For simple logarithmic spirals  the conditions (ii) and (iii)
of the PNSW--hypothesis are not fulfilled because the equation \cite{b1}

\begin{equation}
\label{eq91}
S_{n_1n_2}(x)+S_{n_2n_3}(x)+S_{n_3n_1}(x)=0,
\end{equation}
where
$$
S_{\alpha\beta}(x)=(p_{\alpha} (m) x+1)(p_{\beta}
(m) x+1)\sin\left(\sqrt{\frac{1}{x^2}-1}\;\ln\frac{(p_{\alpha} (m) x+1)}
{(p_{\beta} (m) x+1)}\right)
$$
cannot be satisfied with the same value  $x=\cos\varphi$  for any triplets\linebreak
$(p_{n_1} (m), p_{n_2}(m), p_{n_3}(m))$ from any ray
$r_{m},\, m\in\ov{M}$.

Nevertheless, the solution (\ref{eq91}) for all the denoted prime triplets
from all rays of the  matrix in Appendix 1 shows that $x$ remains
in a sufficiently narrow interval
$I_{x}=(0.202, 0.326)$ with an average
$\ov{x}\approx 0.264$ to which there corresponds a value
$\ov{\varphi}\approx 74.69^{\circ}$.
On Figure 4 a pure--logarithmic web is presented where only the
condition (i) is fulfilled.

\begin{figure}
\centering
\includegraphics[width=9cm, angle=-90]{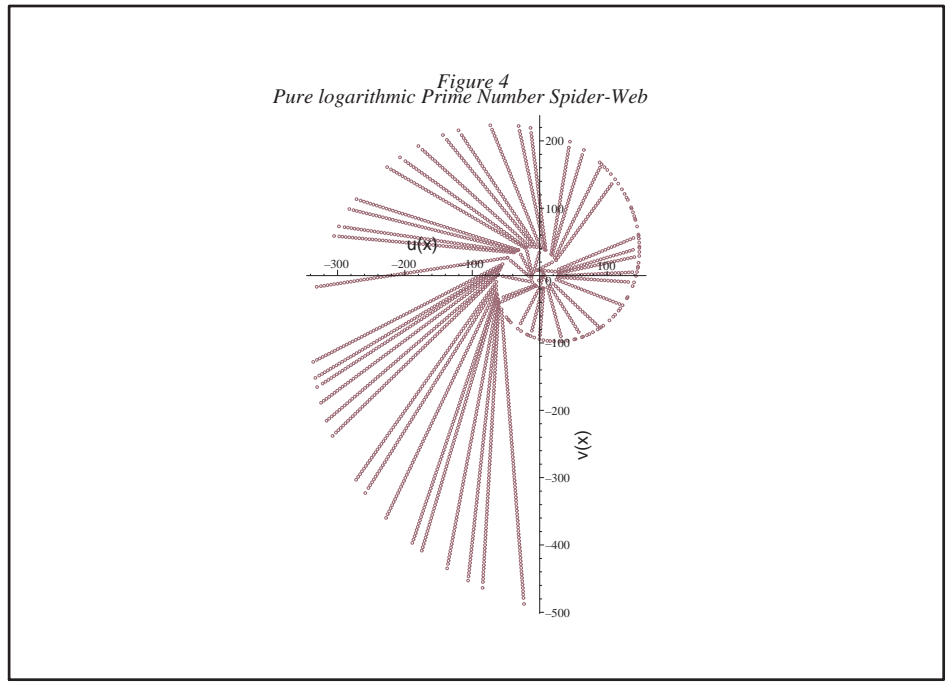}
\caption{{\it Pure logarithmic prime number spider web}}
\end{figure}

This result stimulates us to search for a verification of the
PNSW--hypothesis
in the class of logarithmic spline--spirals (LSS):
$$
\begin{array}{l}
\rho_{s_1}(\theta)=e^{\di s_1(\theta)},\\
\phantom{qqq}\\
s_1(\theta)=\left\{
\begin{array}{l}
\alpha_{i+1}\theta+\beta_{i+1},\quad \theta_i\leq\theta\leq\theta_{i+1}\quad
0\leq i\leq k -1,\\
\alpha_{i-1}\theta_{i-1}+\beta_{i-1}=\alpha_{i}\theta_{i-1}+
\beta_{i},\quad 2\leq i\leq k,\\
\end{array}
\right.
\end{array}
$$
where the first degree spline $s_1 (\theta)$ is defined on an irregular set
$$
\Delta_{k} : 0=\theta_0<\theta_1<\theta_2<\ldots<\theta_{k -1}<
\theta_{k},
$$
with a number $k\geq 3$ of subintervals
$[\theta_i, \theta_{i+1}],\quad 0\leq i\leq k-1$,
which increases with the number of rotations $n$ of the web $W_n(P)$.

The unknowns in the spiral $\rho_{s_1}(\theta)$ are both the knots of the set
$\Delta_k$ and the spline--spiral coefficients of the elements
$e^{\di \alpha_i \theta+\beta_i}$

$$
\{\alpha_i, \theta_i, \beta_i\}_{i=1,2,\ldots, k}.
$$
They are determined from the conditions of the PNSW--hypothesis
with regard for the initial condition
\begin{equation}
\label{eq10}
e^{\di \alpha_1\theta_0+\beta_1}=1\Longrightarrow\beta_1=0.
\end{equation}

{\it For arbitrary  $x\in\mathbb{R}_+^1$ there exists
an unique $p(k_x)\in P$ such that
\linebreak $p(k_x-1)\leq x<p(k_x)$ and the isometric transformation
$x\in\mathbb{R}_+^1$ on $\mathbb{R}^2$, determined by the condition (i)
of the PNSW--hypothesis, acquires the explicit form
\begin{equation}
h_{\di\rho_{s_1}}(x): \mathbb{R}^1_{+}\to \lambda (0, \theta_x)=
p(k_x-1)+
\sqrt{1+\frac{1}{\alpha_{k_x}^2}}\hspace{0.3cm} e^{\di \beta_{k_x}}
\left(  e^{\di \alpha_{k_x}
\theta_x}-e^{\di \alpha_{k_x}\theta_{k_x-1}}\right),
\label{eq11}
\end{equation}
where
$$
\begin{array}{l}
0\leq x<\infty,\quad \theta_x=\frac{\di 1}{\di\alpha_{k_x}}\ln
E(x),\\
\phantom{tttt}\\
E(x)=\frac{\di\alpha_{k_x}}{\sqrt{\di 1+\alpha_{k_x}^2}}
e^{\di -\beta_{k_x}}
(x-p(k_x-1))+e^{\di\alpha_{k_x}\theta_{k_x-1}},\quad p(0)=0.
\end{array}
$$

The points $(u_x, v_x)\in\rho_{s_1}(\theta)\subset\mathbb{R}^2$, which correspond
to the numbers $x\in\mathbb{R}_+^1$, have the Euler and Cartesian
coordinates respectively:

\begin{equation}
\label{eq12}
\rho (\theta_x)=e^{\di \alpha_{k_x}\theta_{x}+\beta_{k_x}}=e^{\di\beta_{k_x}}
\ln E(x),\quad \theta_x=\frac{\di 1}{\di \alpha_{k_x}}
\ln E(x)\\
\end{equation}
\phantom{aaaaaaaaa}and
\begin{equation}
\label{eq13}
u_x=\rho(\theta_x)\cos(\theta_x),\quad v_x=\rho(\theta_x)\sin(\theta_x).
\end{equation}}

The first plane spiral isometric to the semi--axis $\mathbb{R}_+^1$
was constructed in  \cite{b9}.

\subsection{About constructing the webs $W_n$}

Attempts to construct the spiral $\rho_{s_1}(\theta)$ under the
 PNSW-hypothesis for a given $n$ lead to a denial of some number
 $k^0$ of starting primes because of the difficulty in fulfilling the
  condition
(ii) around the origin of $\mathbb{R}^2$ (condition (i)
remains valid for the missed primes). In this paper the case
$k^0=11$ is considered, i.e., instead of the rays
$r_1, r_4, r_6, r_8, r_9$ and $r_{10}$, the truncated rays
$\ov{r}_1, \ov{r}_4,  \ov{r}_8, \ov{r}_9$ and $\ov{r}_{10}$
obey the condition (ii), and these rays start with the numbers
$127, 59, 41, 87, 83$ and $109$ respectively.

According to (\ref{eq10}), to the first element
$e^{\di \alpha_1\theta}\; (0\leq\theta\leq\theta_1)$ of the spiral
$\rho_{s_1}$ there  corresponds the real segment $[0,\, p(k^0+1)]$.

The rotations $W_n$ are taken in account from the ray
$$
r_{12}=\{37, 157, 919, 7193,\ldots\}
$$
in the direction counter--clockwise.

At first, $\rho_{s_1}^{(3)}$  and $W_3$ are constructed on the basis
of the first 3 elements of the first 25 rays
 $\ov{r}_1, \ov{r}_4,\ldots,  r_{36}$ plus the fourth element of the ray
 $r_{12}$.\\

The satisfaction of the conditions (i), (ii) of the   PNSW--hypothesis
with regard to (\ref{eq10}) leads to a solution of a nonlinear system
({\it $W_3$--system}) of 228 equations and 304 inequalities,
caused by the condition (iii), with respect to the 228 unknowns
$$
\{ \alpha_i, \theta_i, \beta_i\}_{i=1,2,\ldots, 76}\;.
$$

 All 832 primes up to the number
7193, which remain unused in the construction of the $W_3$--system
$(p_{-1}(7193)=919;\; 11+ (25\times 3+1)+\linebreak 832=919)$ are placed by
the Cartesian coordinates (\ref{eq13}) on the $2nd$ and
$3rd$ rotation of $\rho_{s_1}^{(3)}$.

\begin{figure}
\centering
\includegraphics[angle=-90, width=16cm, keepaspectratio]{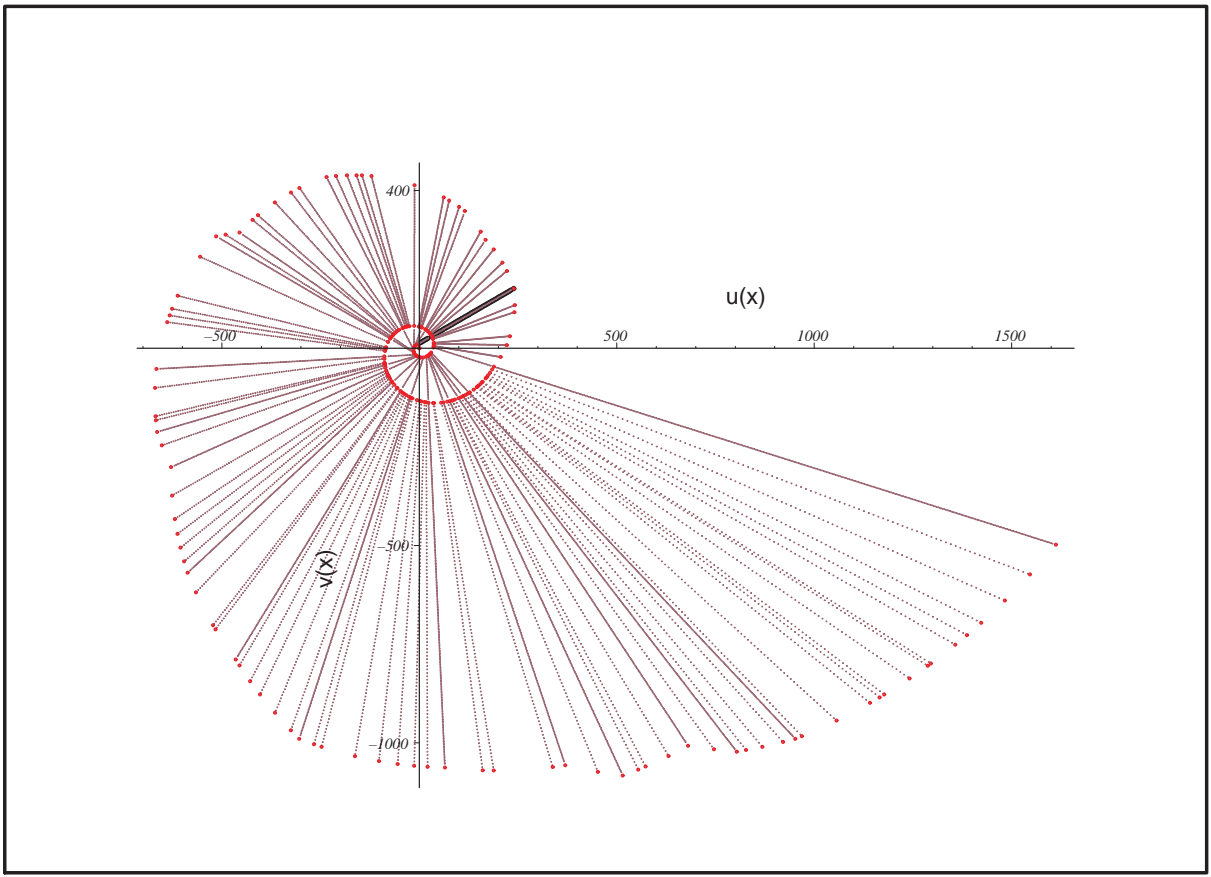}
\caption{{\it 3--rotation prime number spider web.
The thick black line denotes the initial ray $r_{12}$. The direction of
rotation is counter-clockwise}}
\end{figure}

The building up of new rotations $n>3$ on $\rho_{s_1}^{(3)}$ is reduced
to the subsequent solution of  $3\times 3$
nonlinear systems of equations.

Both solvability and uniqueness  of the mentioned infinite set of
nonlinear systems are the analytical interpretation of the content
of the PNSW-hypothesis.\\

A desire to avoid the solution of the  $W_3$--system leads to a
construction of approximations $\widetilde{W}_3, \widetilde{W}_4$ and
$\widehat{W}_3$, in which the first two turns are constructed as a simple spiral
from $\mathscr{L}_{\varphi}$. The subsequent rotations are constructed as  LSS.
In these webs $\varphi=74.18896^{\circ}$.

The web $\widetilde{W}_3$ fairly well illustrates the main properties of the
prime number spider webs and the web $\widetilde{W}_4$ shows
the possibility of continuing the construction of higher rotations.
The web $\widehat{W}_3$ is created for a generation of initial
approximations to a solution of the $W_3$--system.
It also illustrates  the demerits of the approximated webs.

\begin{figure}
\centering
\includegraphics[angle=-90, width=16cm, keepaspectratio]{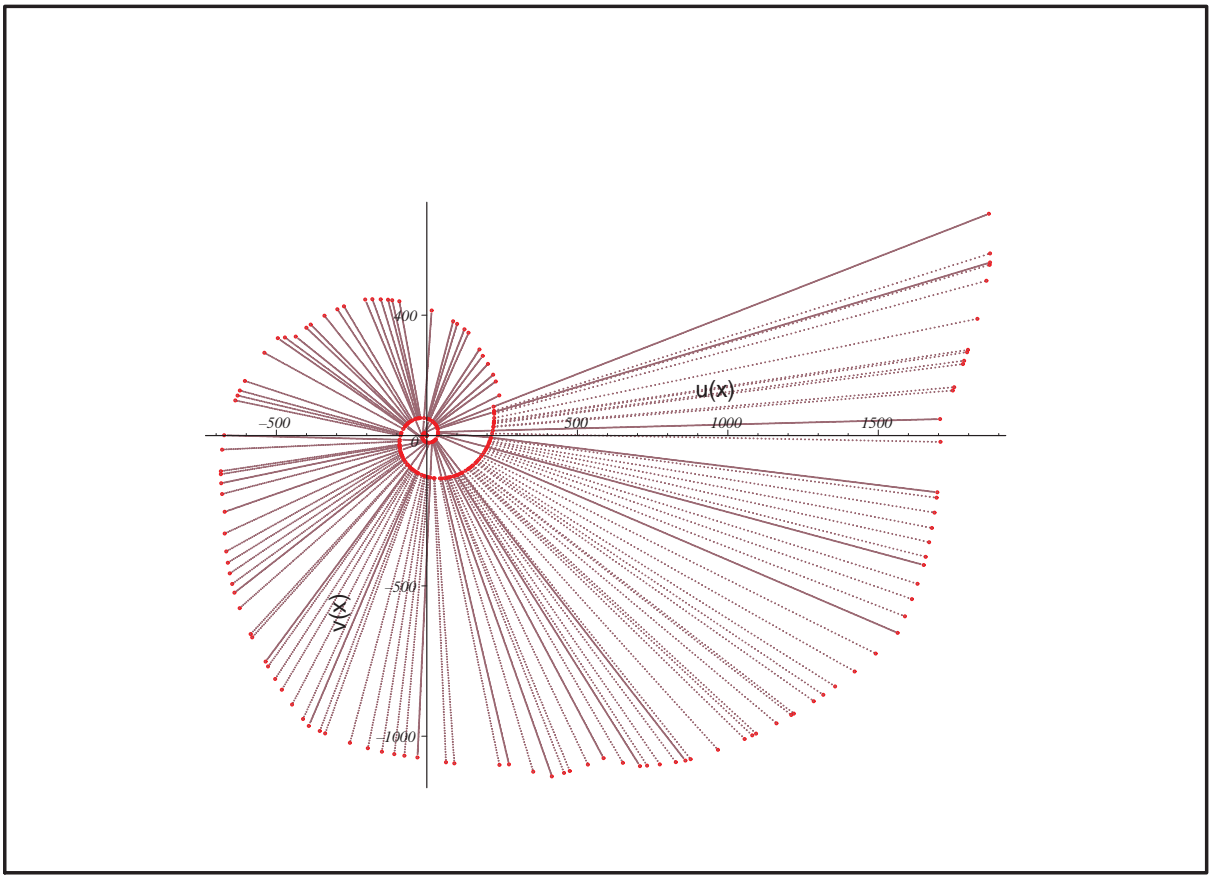}
\caption{{\it 3--rotation prime number spider web with full  $3rd$ turn}}
\end{figure}

The web $\widetilde{W}_3$ is presented in Figure  5. It is constructed on
the basis of 211 primes: 3 numbers from each of the first 20 rays from
$\ov{r}_1$ to $r_{30}$, 2 numbers from 5 rays from $r_{32}$ to $r_{36}$
and 2 numbers from each of the 71 rays from $r_{38}$ to $r_{126}\,
(3\cdot 20+2\cdot 5 +2\cdot 71=211)$; 147 of these primes are placed
on the $2nd$ and  the $3rd$ rotations by means of  coordinates
(\ref{eq13}).

The web $\widetilde{W}_3$ does not have a complete $3rd$ rotation,
ending in the number
5381 from the ray $\ov{r}_1$ and not reaching number 7193 from  the
initial ray $r_{12}$, marked in Figures 5 and 9 by a thick black line.
To simplify the figure, 498 primes remaining until
$p_{-1}(5381)=709$  (709-211=498) are not placed
 on the $3rd$ rotation.

The web $\widehat{W}_3$ with a full $3rd$ rotation is represented in Figure
6. It is built on the basis of 255 primes: 3 from the first 25 rays from
$\ov{r}_1$ to  $r_{36}$ and 2 from 90 rays from  $r_{38}$ to
$r_{151}$\, $(3\cdot 25+2\cdot 90=255)$; 180 of these primes are placed on
 the $2nd$ and the $3rd$ rotations by means of coordinates
(\ref{eq13}).

To simplify  the figure, 664 primes (919-255=664) remaining up to\linebreak
$p_{-1}(7193)=919$  are not placed
on the $3rd$ rotation.

\begin{figure}
\centering
\includegraphics[ angle=-90, width=16cm, keepaspectratio]{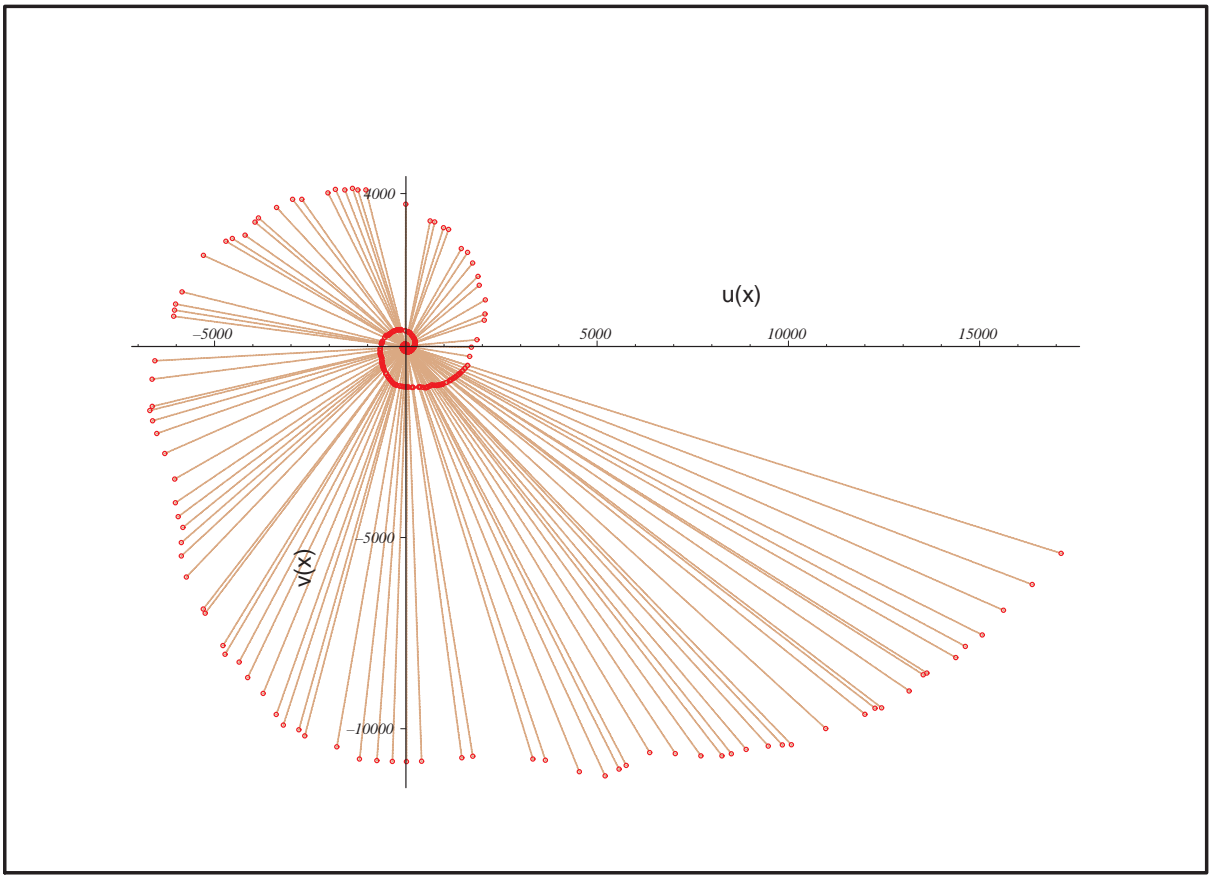}
\caption{{\it 4--rotation prime number spider web. In the center of the picture
one can see the web $\widetilde{W}_3$}}
\end{figure}

At the end of the $3rd$ rotation $\widehat{W}_3$,
(at the transition of the logarithmic spiral in  LSS) a unessential
{\it qualitative defect}
shows up between numbers  877 from $r_{36}$ and 919 from
 $r_{12}$.
The pointed defect
obstructs the exact sewing of the spirals between the numbers
6823 and 7193 from the corresponding rays  $r_{36}$ and $r_{12}$.
This defect is removable  by solving the $W_3$--system.

The web $\widetilde{W}_4$, represented in Figure 7, is obtained from the web
$\widetilde{W}_3$ by adding on the $4th$ rotation as LSS.
In this building, 96 primes are used:
4 elements of the rays from $\ov{r}_1$
to $r_{30}$, 3 elements of the rays  from  $r_{32}$ to $r_{36}$
and 4 elements of the rays from
$r_{38}$ to $r_{126}$ (20+5+71=96).
By means of 116 primes the LSS--units
$e^{\di \alpha_i\theta_i+\beta_i},\quad i=1,2,\ldots,116$
are determined by solving a $3\times 3$ nonlinear system 116
times with respect to
348 unknowns $\alpha_i, \theta_i$ and $\beta_i$:
94 times exactly and  22 times approximately.
The cases of inexact solutions
result in unessential distortions of the condition (ii) for   22 rays
(in the diagrams of  $\widetilde{W}_3$ and $\widetilde{W}_4$ those distortions
are not seen).

\begin{figure}
\centering
\includegraphics[width=13cm, totalheight=11cm]{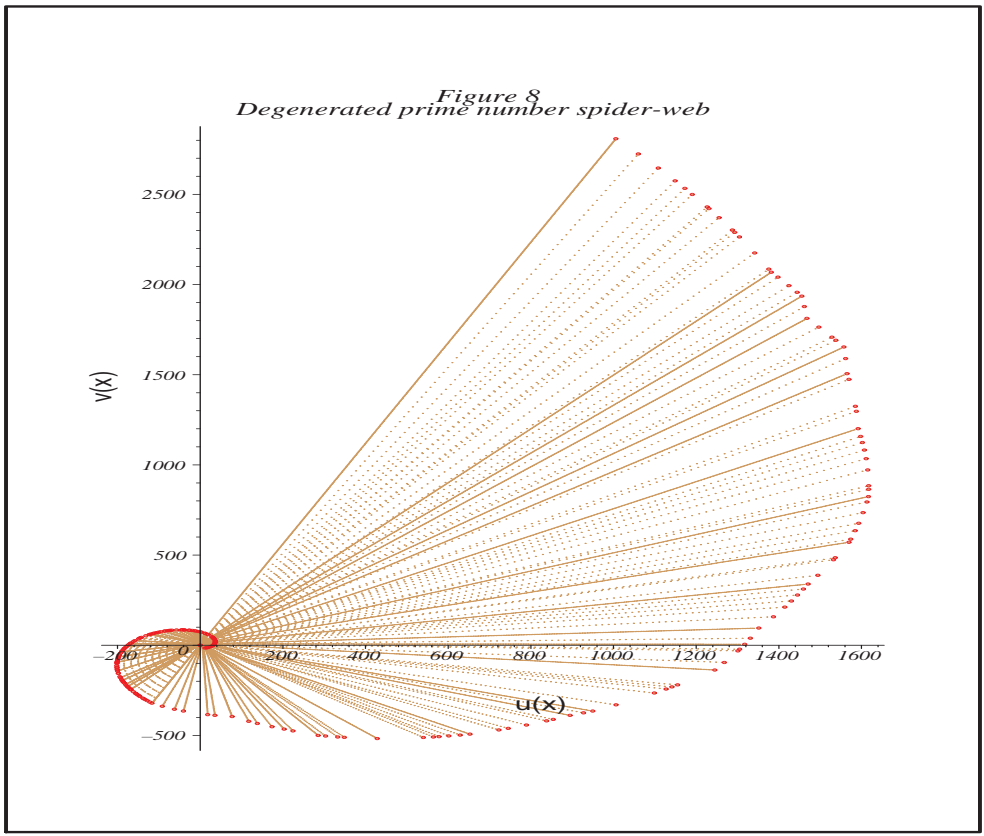}
\caption{{\it Degenerated prime number spider web}}
\end{figure}

A variety of possibilities for constructing the  PNSW--hypothesis is
illustrated by web $\widetilde{W}_{deg}$ (Figure 8),
in which the condition (ii) is violated in the following way:
the rays lie on straight lines, but the subsequent segments of the rays have
the opposite directions.
In $\widetilde{W}_{deg}$, the same primes are used as in
$\widetilde{W}_{3}$.\\

\subsection{Web formation rules and properties}

Satisfying the PNSW--hypothesis, the conditions of  $\widetilde{W}_{3}$ and
$\widetilde{W}_{4}$ fix the individual peculiarities of the behaviour
of primes around the origin of $\mathbb{R}^2$ as concrete systems of embedded
trapezoids confined between the rays
$\ell_{p_{\mu 1}}$ and $\ell_{p_{(\mu+k)1}},\; k\geq 1$.
For example, the trapezoid $t_{19}$ confined between the rays
$\ell_{113}$ and $\ell_{127}$ (Figure 9) is one of the nineteen
{\it 3-rotation embedded trapezoids (3RET)} of the web $\widetilde{W}_{3}$.\\

From the Theorem 2 it  follows that

\begin{itemize}
 \item[$\mathbf{q_8)}$]   The elements of the clusters
 $c_{\alpha_{\mu}}(\alpha_{\mu})$   (they and only they)
become the starts of new rays on the spiral $\rho_{s_1}(\theta)$.\\

Let a cluster
$$
c_{\mu}(k)=\{p_{\mu 1}, p_{(\mu+1)1},\ldots, p_{(\mu+k)1}\}.
$$
lie on
the $\nu$th rotation of the spiral  $\rho_{s_1}^{(n)}(\theta)$.
On the  $(\nu+q)$th rotation to it there  correspond the primes
$\{p_q(p_{(\mu+i)1})\}_{i=0,1,\ldots, k}$.

The PNSW--hypothesis shows that

\item[$\mathbf{q_9)}$]   The geometric figure on $\mathbb{R}^2$,
confined between the arcs\linebreak $(p_{\mu 1},p_{(\mu+k)1})$ and
$(p_q(p_{\mu 1}),p_q(p_{(\mu+k)1}))$ by
the spiral $\rho_{s_1}$ and the segments $|p_{\mu 1},p_q(p_{\mu 1})|$,
$|p_{(\mu +k)1},p_q(p_{\mu+k})|$ of the rays $\ell_{p_{\mu 1}}$ and
$\ell_{p_{(\mu+k)1}}$, is the concave--convex trapezoid
$$
z(\nu, \mu, k, q)=[(p_{\mu 1},p_{(\mu +k)1}),\; (p_q(p_{\mu 1}),
p_q(p_{(\mu+k)1}))].
$$

The trapezoids of the type $z(\nu, \mu, 1,1)$ are {\it elementary
trapezoids (or holes)} to $W_n$.
For example
$$
z(1, 19, 1, 1)=[(113, 127),\, (617, 709)]
$$
is an elementary $W_2$--trapezoid (Figure 10).

\begin{figure}
\centering
\includegraphics[angle=-90, width=14cm, keepaspectratio]{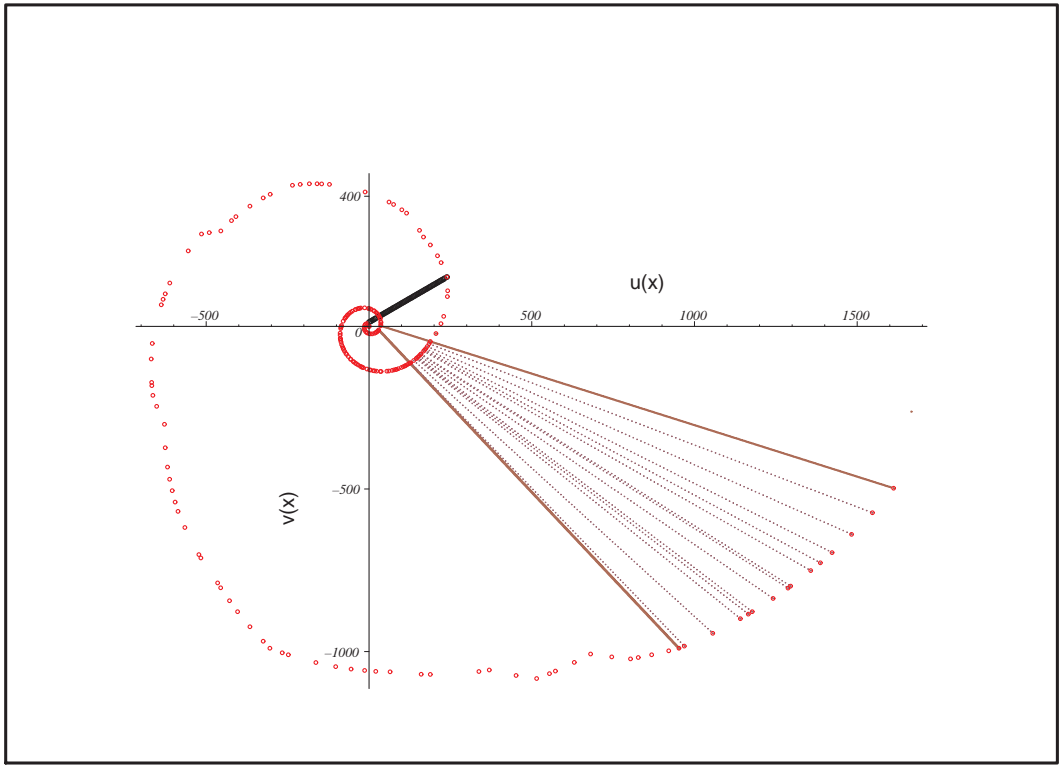}
\caption{{\it 3--rotation system of embedded trapezoids
$3RET_{19}=z(1,19,1,2)=[(113, 127),\, (4549, 5381)]$.
 The thick black line
defines the initial ray $r_{12}$. The direction of rotation is
counter-clockwise}}
\end{figure}

\begin{figure}
\centering
\includegraphics[angle=-90, width=11.5cm, keepaspectratio]{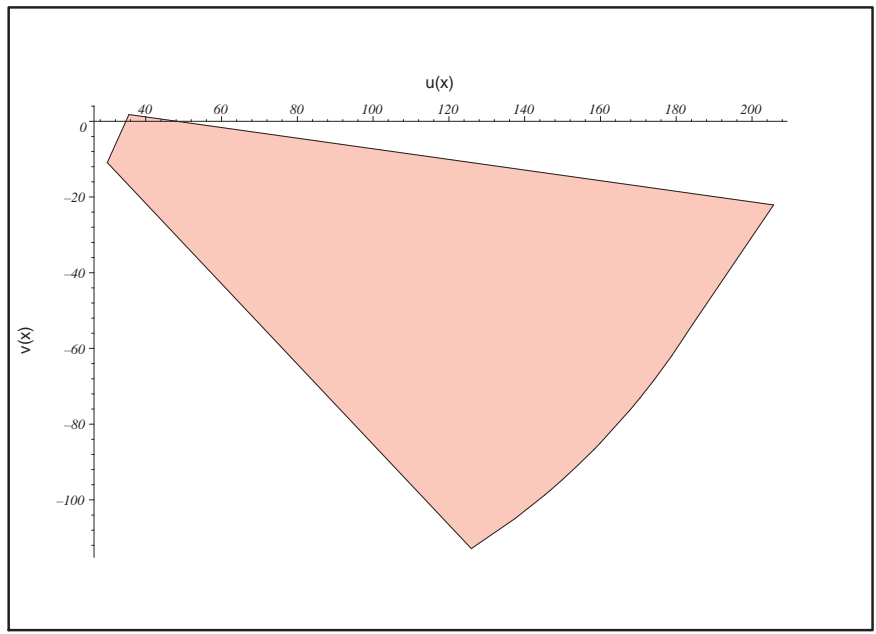}
\caption{{\it The elementary trapezoid
$z_{19}=z(1,19,1,1)=[(113, 127),\, (617, 709)]$}}
\end{figure}

Let  the primes
 $p_{\mu 1},$ $p_{(\mu+k)1}$ $\in c_{\mu}(\alpha_{\mu}),
\alpha_{\mu}\geq 3,\; 0\leq k\leq\alpha_{\mu}$ lie
on the $\nu$th rotation of the spiral $\rho_{s_1}^{(n)}$.
On the $(\nu+1)$th rotation,
to them there corresponds the cluster $c_{\mu_1}(\alpha_{\mu_1})=
\{p_{\mu_1 1}, p_{(\mu_1+1)1},\ldots,
p_{(\mu_1+\alpha_{\mu_1})1}\}$ with a length $\alpha_{\mu_1}=
p_{(\mu+k)1}-p_{\mu 1}-1$, according to Corollary 3.

From the conditions of the PNSW-hypothesis there stems
the following rule for formation of 3RET:

 \item[$\mathbf{q_{10})}$]  On the plane $\mathbb{R}^2$ the following
 equalities hold

\begin{equation}
\label{eq23}
 z(\nu+1,\mu_1, \alpha_{\mu_1}+1, 1)=\bigcup
 \limits_{i=0}^{\alpha_{\mu_1}}z(\nu+1,\mu +i, 1, 1),
\end{equation}
\begin{equation}
\label{eq24}
 z(\nu,\mu, 1, 2)=z(\nu,\mu, 1, 1)\cup z(\nu+1, \mu_1,
 \alpha_{\mu_1}+1, 1).
\end{equation}

Equalities (\ref{eq23}) and (\ref{eq24}) can be regarded as those
between the areas (the sign $\cup$ changes to  +).

Let us apply the rule for formation of 3RET to $3RET_{19}$.
Using the matrix from Appendix 1, we obtain
$$
3RET_{19}=[(113, 127),\, (617, 709)]\cup [(617, 619),\, (4549, 4567)]\cup
$$
$$
[(619, 631),\, (4567, 4663)]\cup\ldots \cup
 [(701, 709),\, (5281, 5381)].
$$
In this example $\nu=1,\, \mu=30, \, \alpha_{\mu_1}=127-113-1=13$.
The starts of the newly--appeared rays on the second rotation are
$$
\begin{array}{llll}
p(114)=619,& p(115)=631,& p(116)=641,& p(117)=643,\\
p(118)=647,& p(119)=653,& p(120)=659,& p(121)=661,\\
p(122)=673,& p(123)=677,& p(124)=683,& p(125)=691, \mbox{ and }\\
p(126)=701.\\
\end{array}
$$

\item[$\mathbf{q_{11})}$] The rule for formation of 3RET shows that 3RET
is a mosaic of $\alpha_{\mu_1}+1$ elementary trapezoids, which
are grouped in the direction of the ray $\ell_{p_{\mu 1}}$
in the following way:

the first is the elementary trapezoid
$$
[(p_{\mu 1},\, p_{(\mu+1)1}),\,(p_1(p_{\mu 1}),\, p_1(p_{(\mu+1)1}))],
$$
 followed by the composite trapezoid
$$
[(p_1(p_{\mu 1}),\, p_1(p_{(\mu+\alpha_{\mu_1})1})),\,
(p_2(p_{\mu 1}),\, p_2(p_{(\mu+\alpha_{\mu_1})1}))].
$$

The elements (elementary trapezoids too) of the last one cause in total
$\mu_2= p_1(p_{(\mu+\alpha_{\mu_1})1})-p_1(p_{\mu 1})-1$ new rays
on the $(\nu+2)th$ rotation,
to which again the 3RET formation rule is applied for obtaining
$\alpha_{\mu_2}$  new composite trapezoids located between the
$3rd$ and $4th$ rotations of $W_n$.

Under the condition of the PNSW-hypothesis the process of growth
of the initial 3RET is unlimited and leads to the formation of {\it a class
$(\ell_{p_{\mu 1}},\,\ell_{p_{(\mu+1)1}})$--trapezoids}
containing an unlimited number of different elementary trapezoids.

\item[$\mathbf{q_{12})}$] Let  {\it mitos}   denote a closed geometric
figure locked between the arc\linebreak $(p(12), p_{1}(12)))$ of the spiral
$\rho_{s_{1}}$ and the segment $|p(12), p_1(12)|$ of the ray $r_{12}$.
Then the plane $\mathbb{R}^2$ can be represented as a mosaic of
the initial $k(k^0)$--classes of embedded trapezoids
(k(11)=25):

$$
\mathbb{R}^2=\bigcup\limits_{i=1}^{k(k^0)}\left(
\ell_{p(k^0+i)},\, \ell_{p(k^0+i+1)}\right)\cup \mbox{{\it mitos}}.
$$

Inside the  {\it mitos}  the starts of the first few rays can intersect
each other, i.e., be in a mitosis status.

Properties $q_{10}),\, q_{11})$ and $q_{12})$  are generalized
in the property $q_{13}$.

\item[$\mathbf{q_{13})}$]  The plane $\mathbb{R}^2$ is representable
as a mosaic  (in the theoretical--set sum sense)
of all the prime--numerical elementary trapezoids and {\it mitos}.

\item[$\mathbf{q_{14})}$]  The PNSW-hypothesis leads to the geometric
interpretation of $\pi (x)$:

\begin{equation}
\label{eq16}
\pi (p_{\nu}(m))=p_{\nu-1}(m)
=\lambda (0,\, p_{\nu-1}(m)),\quad \nu\geq 2,\, m\in \ov{M}.
\end{equation}

Equalities (\ref{eq16}) allow one to interpret  geometrically
the  $\Omega$--theorem of Littlwood  and Theorem 1 from \cite{b10}, where
for setting the real values and in particular the values of
$Li(x)$ again the transformation $h_{\di \rho_{s_1}}(x)$
and coordinates (\ref{eq13}) are used.

The geometric interpretation of the Theorem 3 consists
in determination of the laws for  appearance of
 $u$-- and $b$--twins on $W_n$.

\item[$\mathbf{q_{15})}$]  $u$--twins  appear on $W_n$ in the following
cases:

 a)  as starts of new rays $\ell_{p_{\mu 1}}$ and $\ell_{p_{(\mu+1)1}}$
on  the $1st$ rotation  3RET, which on the $2nd$ rotation cause one new ray
$\ell_{p((p_{\mu 1}+p_{(\mu+1)1})/2)}$. Then  3RET  is composed of three elementary trapezoids.

Examples:
$$
\ov{m}_8(3)\to(\ell_{71}, \ell_{73})\to\ell_{359},
$$
$$
\ov{m}_9(5)\to(\ell_{101}, \ell_{103})\to\ell_{557},
$$
$$
\ov{m}_{11}(5)\to\left\{
\begin{array}{l}
(\ell_{137}, \ell_{139})\to\ell_{787},\\
(\ell_{149}, \ell_{151})\to\ell_{863};
\end{array}
\right.
$$

 b)  as a pair of subsequent primes on the $2nd$ rotation of 3RET

Examples on the $2nd$ turn of the trapezoid $z(1, 19, 1, 2)$ (Figure 9):

$$
\ov{m}_{30}(13)\to\left\{
\begin{array}{l}
(\ell_{641}, \ell_{643})\to\ell_{4783},\\
(\ell_{659}, \ell_{661})\to\ell_{4937}.
\end{array}
\right.
$$

\item[$\mathbf{q_{16})}$]   One of the elements  of $b$--twin
$(t_1(\ov{n}),\, t_2(\ov{n}))$
lies on the existing ray\linebreak  $\ell_{p_q(p_{\mu_1 1})},\; q\geq 1$,
but the other element is a start of a new ray $\ell_{p_{\mu_2 1}}$:
$$
t_1(\ov{n})\equiv p_q(p_{\mu_1 1}),\, t_2(\ov{n})\equiv p_{\mu_2 1}
\mbox{ -- $right$ twin};
$$
$$
t_1(\ov{n})\equiv p_{\mu_2 1},\, t_2(\ov{n})\equiv p_q(p_{\mu_1 1})
\mbox{ -- $left$ twin}
$$

 In such a way the $b$--elementary trapezoid
$$
[(p_q(p_{\mu_1 1}), p_{\mu_2 1}),\; (p_{q+1}(p_{\mu_1 1}), p(p_{\mu_2 1}))]
$$
 is sewn together
to the right of the ray $\ell_{p_{\mu_1 1}}$, and trapezoid
$$
[(p_{\mu_2 1},p_q(p_{\mu_1 1})),\; (p(p_{\mu_2 1}),p_{q+1}(p_{\mu_1 1}))]
$$
 to the left of the ray  $\ell_{p_{\mu_1 1}}$.

Examples: the trapezoid
$[(617, 619),\,(4549, 4567)]$ is sewn to the right of the ray
$\ell_{113}\,\left((\ell_{617}, \ell_{619})\to\ell_{4561}\right);$

the trapezoid
$[(857, 859),\,(6653, 6661)]$ is sewn to the left of the ray\linebreak
$\ell_{149}\, \left((\ell_{857}, \ell_{859})\to\ell_{6659}\right)$.

The properties $q_{15})$ and  $q_{16})$ show the following $W_n$ sewing
property

\item[$\mathbf{q_{17})}$]  The twin pairs $sew$ uniformly  over $n$
the Eratosthenes rays
in an $unified$ plane web  $W_n,$ $n=1,2,\ldots\; .$
\end{itemize}

\section{Conclusion}

The study of  the inner prime number distribution law remains at the
initial stage.

Finally, we shall note how the proof of Conjecture 1 looks like,
and we shall find a  possibility of its generalization.
We should like to indicate possible
applications of the proposed, in this paper, approach to a study of
the oddish prime number behaviour.\\

At the first stage of the proof of Conjecture 1
it is necessary to solve the
$W_3$--system and  revise the number $k^0$.
If the $W_3$--system cannot be solved with $k^0=11$, then
the numbers $k^0=10$ and $k^0=12$ should be tried.

In the formulation of the $W_3$--system a function
$\tilde{p}(x)\in C^1 (0,10000],$
  with a property
\begin{equation}
\label{eq20}
\tilde{p}(n)= p(n),\quad n=1,2,\ldots,1229
\end{equation}
is used.

This function can be built up by  {\it the rod spline method}
\cite{b11} using as a rod the approximation (\cite{b12}, exercise 9.21)
$$
\ov{p}(x)=x\left(\ln x +\ln\ln x+\frac{\ln\ln x-2}{\ln x}-
\frac{{(\ln\ln x)}^2/2-3\ln\ln x+5.5}{{(\ln x)}^2}-1\right).
$$

The searched function will be of the form  $\tilde{p}(x)=s_2(x)\ov{p}(x)$,
where $s_2(x)$ is a parabolic spline determined on two nonuniform sets.
At the interpolation points equalities (\ref{eq20}) are fulfilled.

The  $W_3$--system solvability  can first
be investigated  numerically
 by using, for example, the program LANCELOT \cite{b13}.

At the second stage of the proof of Conjecture 1, it is necessary to prove
inductively the continuability  of the {\it basic $k(k^0)$--classes of
embedded trapezoids}
$$
(\ell_{p_{(k^0+i)1}},\,\ell_{p_{(k^0+i+1)1}}),\quad i=1,2,\ldots, k(k^0),
$$
as well as the continuability of {\it the new classes on rotations} $n\geq 2$.\\

Conjecture 1 can be extended to all {\it mesm$_f$}--matrices.\\


{\bf Conjecture 2.}
{\it For each matrix $A_f\in\mathscr{E}_f$ there exists
 LSS-spiral  and a web $W_n (A_f)$, satisfying the conditions
(i), (ii) and (iii) of the  PNSW--hypothesis. }\\

In the infinite set of webs $\mathbf{w}=
\{W_n(A_f):A_f\in\mathscr{E}_f,\; n=1,2,\ldots\}$
it is necessary to introduce an operation  {\it summing webs} @
in such a way  that the equalities
$$
W_n(S)@W_n(T_1)@W_n(T_2)=W_n(P)
$$
and
$$
W_n(D_{6n-1})@W_n(D_{6n+1})=W_n(P),
$$
hold, by analogy with the set-theoretical  equalities
$P=S\cup T_1\cup T_2$ and
 $P=D_{6n-1}\cup D_{6n+1}\cup \{2, 3\},\; \mbox{ where}
 \{2, 3\}\subset \mbox{ {\it mitos}}$.

On the whole, the algebraic structures of mesm-matrices
$[B_f\; ^2A_f]$, webs $W_n(A_f)$ and the set
$\mathbf{w}$ remain unexplored.\\

The spiral  $\rho_{s_1}(\theta)$ length  from the
PNSW--hypothesis
can turn out to be an appropriate {\it time coordinate} in the description
of physical processes taking place in  asymmetric and irreversible
time.

Indeed, the mapping (\ref{eq11}) can be extended also for negative values
\linebreak
$x\in (-\infty, 0]:$
$$
h_{\di\rho_{s_1}}(x): \mathbb{R}_{-}^1\to\lambda (0, \theta_x)=
\frac{1}{\cos\varphi}\left(e^{\di (\ctg\varphi)\theta_x}-1\right),
$$
$$
-\infty <x\leq 0,\, \varphi =\arccot (\alpha_1).
$$

In such a way, in the unit circle (inside the domain {\it mitos})
there remains a finite
{\it negative moustache} with a length
$$
\lambda (0, -\infty)=-\frac{1}{\cos\varphi}.
$$

Now, the arc length of the spiral
$\rho_{s_1}(\theta),\, -\infty<\theta<\infty$
is split up in three pieces: the length of {\it the negative
moustache}, the finite positive length spiral in the {\it mitosis}
and the infinite arc length corresponding to the semiaxis
$[p(k^0+1),\infty)$.\\

The solution of equations  (\ref{eq88}), (\ref{eq99})
provides a motivation for the following hypothesis:\\

{\bf Conjecture 3.}
{\it For arbitrary $m\in\ov{M}$ and large $n$ the inequality
\begin{equation}
\left|L(p_{n+1}(m))-p_n(m)\right|\leq
c_1\sqrt{p_{n+1}(m)}\ln (p_{n+1}(m))
\label{eq118}
\end{equation}
is fulfilled with a constant $c_1$  independent of $m$ .}

Inequalities (\ref{eq118}) result in the common estimate
\begin{equation}
\label{last}
\left| L(x)-
\pi (x)\right|<c_2\sqrt{x}\ln x,\; c_2=const.
\end{equation}

The proof of inequality (\ref{eq118}) is easier than the proof
of inequality (\ref{last}). If from estimation (25) follows the truth
of the Riemann hypothesis that complex solutions
of the equations $\zeta (s)=0,\; s\in \not\hspace{-0.15cm}C$ have a form
$s=1/2+i\gamma_n,\;$ $\gamma_n\in\mathbb{R},\; n=1,2,\ldots,$
then the inner prime number distribution law
(9)  will prove to be a new useful
tool of the analytical number theory.

\newpage

{\normalsize
{\bf Appendix 1}\\

{\bf Matrix $\mathbf{[\ov{M}\phantom{d} ^2P]}$}
$$
\begin{array}{llllllllllll}
\hline\\
  1&  2&  3& 5& 11& 31\\
 \, & \, &  127&  709& 5381&  52711\\
\, & \, &  648391& 9737333 & 174440041& 3657500101\\
\, & \, & 88362852307& 2428095424619\\
\,& \,&  75063692618249...\\
 4& 7& 17&  59&  277& 1787\\
 \, & \, & 15299& 167449& 2269733& 37139213\\
 \, & \, &  718064159& 16123689073&  414507281407\\
 \, & \, & 12055296811267...\\
6&   13&  41& 179& 1063& 8527\\
\, & \, & 87803& 1128889& 17624813& 326851121\\
\, & \, & 7069067389& 175650481151&  4952019383323...\\
8&   19&  67& 331& 2221& 19577\\
\, & \, & 219613& 3042161& 50728129& 997525853\\
\, & \, & 22742734291  & 592821132889&  17461204521323...\\
9&   23&  83& 431& 3001& 27457\\
\, & \, & 318211& 4535189& 77557187&  1559861749\\
\, & \, & 36294260117& 963726515729& 28871271685163...\\
10&   29& 109& 599& 4397& 42043\\
\, & \, & 506683& 7474967& 131807699& 2724711961\\
\, & \, & 64988430769& 1765037224331& 53982894593057...\\
12&   37& 157& 919& 7193& 72727\\
\, & \, & 919913& 14161729& 259336153&  5545806481\\
\, & \, & 136395369829& 3809491708961...\\
14&   43& 191& 1153& 9319& 96797\\
\, & \, & 1254739& 19734581& 368345293& 8012791231\\
\, & \, & 200147986693& 5669795882633...\\
15&   47& 211& 1297& 10631& 112129\\
\, & \, & 1471343& 23391799&  440817757& 9672485827\\
\, & \, & 243504973489& 6947574946087...\\
16&   53& 241& 1523& 12763& 137077\\
\, & \, & 1828669& 29499439& 563167303& 12501968177\\
\, & \, & 318083817907& 9163611272327...\\
18&   61& 283& 1847& 15823& 173867\\
\, & \, & 2364361& 38790341& 751783477& 16917026909\\
\, & \, & 435748987787& 12695664159413...\\
20&   71& 353& 2381& 21179& 239489\\
\, & \, & 3338989& 56011909& 1107276647& 25366202179\\
\, & \, & 664090238153& 19638537755027...\\
21&   73& 367& 2477& 22093& 250751\\
\, & \, & 3509299& 59053067&  1170710369& 26887732891\\
\, & \, & 705555301183& 20909033866927...\\
\end{array}
$$

{\bf Appendix 1}\\

{\bf Matrix $\mathbf{[\ov{M}\phantom{d} ^2P]}$}\ldots\quad Continuation 1.
$$
\begin{array}{llllll}
\hline\\
22&   79& 401& 2749& 24859& 285191\\
\, & \, & 4030889& 68425619& 1367161723& 31621854169\\
\, & \, & 835122557939& 24894639811901...\\
24&   89& 461& 3259& 30133& 352007\\
\, & \, & 5054303& 87019979& 1760768239& 41192432219\\
\, & \, & 1099216100167& 33080040753131...\\
25&   97& 509& 3637& 33967& 401519\\
\, & \, & 5823667& 101146501& 2062666783& 48596930311\\
\, & \, & 1305164025929& 39510004035659...\\
26&  101& 547& 3943& 37217& 443419\\
\, & \, & 6478961& 113256643& 2323114841& 55022031709\\
\, & \, & 1484830174901& 45147154715447...\\
27&  103& 563& 4091& 38833& 464939\\
\, & \, & 6816631& 119535373& 2458721501& 58379844161\\
\, & \, & 1579041544637& 48112275898789...\\
28&  107& 587& 4273& 40819& 490643\\
\, & \, & 7220981& 127065427& 2621760397& 62427213623\\
\, & \, & 1692866362237& 51702420222709...\\
30&  113& 617& 4549& 43651& 527623\\
\, & \, & 7807321& 138034009&  2860139341& 68363711327\\
\, & \, & 1860306318433& 56997887937671...\\
  32&  131& 739& 5623& 55351& 683873\\
  \, & \, & 10311439& 185350441& 3898093877& 94434956839\\
  \, & \, & 2606906998739& 80783250929599...\\
    33&  137& 773& 5869& 57943& 718807\\
  \, & \, &   10875143& 196100297&  4135824247& 100450108949\\
  \, & \, & 2773622459039& 86127342906779...\\
    34&  139& 797& 6113& 60647& 755387\\
    \, & \, & 11469013& 207460717&  4387715993& 106839327589\\
    \, & \, & 2956887579073...\\
    35&  149& 859& 6661& 66851& 839483\\
    \, & \, & 12838937& 233784751& 4973864561& 121763369327\\
    \, & \, & 3386468161121...\\
 36&  151& 877& 6823& 68639& 864013\\
 \, & \, & 13243033& 241568891&  5147813641& 126206581463\\
 \, & \, & 3514741569337...\\
 38&  163& 967& 7607& 77431& 985151\\
 \, & \, & 15239333& 280256489& 6016014239& 148471002899\\
 \, & \, & 4159843299587...\\
\end{array}
$$

\newpage

{\bf Appendix 1}\\

{\bf Matrix $\mathbf{[\ov{M}\phantom{d}^2P]}$}\ldots\quad Continuation 2.
$$
\begin{array}{llllll}
\hline\\
 39&  167& 991& 7841& 80071& 1021271\\
 \, & \, & 15837299& 291905681& 6278569691& 155231019913\\
 \, & \, & 4356423418499...\\
40&  173& 1031& 8221& 84347& 1080923\\
\, & \, & 16827557& 311234591& 6715304579& 166500464477\\
\, & \, & 4684808232443...\\
42&  181& 1087& 8719& 90023& 1159901\\
\, & \, & 18143603& 337033877& 7300206493& 181639026043\\
\, & \, & 5127173445557...\\
    44&  193& 1171& 9461& 98519& 1278779\\
    \, & \, & 20137253& 376292689&  8194134017& 204869160779\\
    \, & \, & 5808496248769...\\
   45&  197& 1201& 9739& 101701& 1323503\\
   \, & \, & 20890789& 391182829& 8534307629& 213736527847\\
   \, & \, & 6069307408303...\\
  46&  199& 1217& 9859& 103069& 1342907\\
  \, & \, & 21219089& 397681327&  8682977119& 217616274683\\
  \, & \, & 6183541562551...\\
   48&  223& 1409& 11743& 125113& 1656649\\
   \, & \, & 26548261& 503859997& 11126538823& 281736685679\\
   \, & \, & 8081022964981...\\
    49&  227& 1433&  11953&  127643&  1693031\\
    \, & \, &  27170047& 516340703&  11415461989&  289357897711\\
    \, & \, &  8307635814431...\\50&  229& 1447& 12097& 129229& 1715761\\
    \, & \, & 27560453& 524172379& 11596829689& 294145810687\\
    \, & \, & 8450108859131...\\
    51&  233& 1471& 12301& 131707& 1751411\\
    \, & \, & 28171007& 536433767&  11881126321& 301656862553\\
    \, & \, & 8673774992821...\\
    52&  239& 1499& 12547& 134597& 1793237\\
    \, & \, & 28889363& 550881943&  12216514841& 310526940547\\
    \, & \, & 8938160481557...\\
    54&  251& 1597& 13469& 145547& 1950629\\
    \, & \, & 31599859& 605555557& 13489097663& 344268078839\\
    \, & \, & 9946200971687...\\
    55&  257& 1621& 13709& 148439& 1993039\\
    \, & \, & 32332763& 620393003&  13835380799& 353471438263\\
    \, & \, & 10221768670013...\\
    56&  263& 1669& 14177& 153877& 2071583\\
    \, & \, & 33691309& 647927381&  14478972721& 370600719481\\
    \, & \, & 10735307868743...\\
\end{array}
$$

{\bf Appendix 1}\\

{\bf Matrix $\mathbf{[\ov{M}\phantom{d} ^2P]}$}\ldots\quad Continuation 3.
$$
\begin{array}{llllll}
\hline\\
 57&  269& 1723& 14723& 160483& 2167937\\
    \, & \, & 35368547& 682005953&  15277169617& 391886115431\\
    \, & \, & 11374585999793...\\
   58&  271& 1741& 14867& 162257&  2193689\\
    \, & \, & 35815873& 691097513& 15490445177& 397580778799\\
    \, & \, & 11545824668459...\\
 60&  281& 1823& 15641& 171697& 2332537\\
  \, & \, &  38235377& 740436923&  16649917331& 428592846379\\
   \, & \, &  12479807093519...\\
   62&  293& 1913& 16519& 182261& 2487943\\
    \, & \, &     40951019& 796000427& 17959785803& 463728180431\\
     \, & \, &  13540770614753...\\
    63&  307& 2027& 17627& 195677& 2685911\\
     \, & \, & 44432569& 867503173& 19651365719& 509248998611\\
  \, & \, &     14919411840803...\\
    64&  311& 2063& 17987& 200017& 2750357\\
     \, & \, &  45564719& 890830471& 20204583739& 524169678691\\
 \, & \, &       15372235794151...\\
    65&  313& 2081& 18149& 202001& 2779781\\
     \, & \, &  46082987& 901517753&  20458245581& 531016168117\\
 \, & \, &       15580165580489...\\
    66&  317& 2099& 18311& 204067& 2810191\\
     \, & \, &  46620709& 912598217&  20721384791& 538121923037\\
      \, & \, &  15796066509169...\\
    68&  337&  2269& 20063& 225503& 3129913\\
     \, & \, &  52286593&  1029838717& 23513901553& 613739626127\\
      \, & \, &  18099406558319...\\
       69&  347&  2341& 20773& 234293& 3260657\\
       \, & \, & 54615469& 1078227191& 24670634249& 645165616243\\
       \, & \, & 19059563752283...\\
   70&  349&  2351& 20899& 235891& 3284657\\
   \, & \, & 55043683&  1087126459& 24883634693& 650958710863\\
   \, & \, & 19236734782351...\\
   72&  359&  2417& 21529& 243781& 3403457\\
   \, & \, & 57160969& 1131224411& 25940205719& 679722101701\\
   \, & \, & 20117195040149...\\
 74&  373& 2549& 22811& 259657& 3643579\\
 \, & \, & 61460533&  1221036307& 28097383163& 738585245417\\
 \, & \, & 21922891272739...\\
  75&  379&  2609& 23431& 267439& 3760921\\
  \, & \, & 63567289& 1265161649& 29159843309& 767640499331\\
  \, & \, & 22816010162129...\\
\end{array}
$$

{\bf Appendix 1}\\

{\bf Matrix $\mathbf{[\ov{M}\phantom{d} ^2P]}$}\ldots\quad Continuation 4.
$$
\begin{array}{llllll}
\hline\\
  76&  383& 2647& 23801& 271939&    3829223\\
   \, & \, & 64795981&  1290918281& 29780778613& 784640376427\\
   \, & \, & 23339094889519...\\
    77&  389& 2683& 24107& 275837& 3888551\\
    \, & \, & 65864459&  1313343397& 30321784529& 799462887341\\
    \, & \, & 23795492951147...\\
 78&  397& 2719& 24509& 280913& 3965483\\
    \, & \, & 67247771&  1342401539& 31023447269& 818701472243\\
    \, & \, & 24388288001989...\\
    80&  409& 2803& 25423& 292489& 4142053\\
    \, & \, & 70432519&  1409422013& 32644249103& 863205467819\\
    \, & \, & 25761357737977...\\
    81&  419&  2897& 26371& 304553& 4326473\\
    \, & \, & 73768631&  1479780677& 34349423377& 910115902141\\
    \, & \, & 27211243680073...\\
    82&  421&   2909& 26489& 305999& 4348681\\
    \, & \, & 74172503& 1488302867& 34556157661& 915809403721\\
    \, & \, & 27387388206553...\\
   84&  433& 3019& 27689& 321017& 4578163\\
   \, & \, & 78339559&  1576442723& 36697520357& 974856473813\\
   \, & \, & 29216297536511...\\
85&  439&  3067& 28109& 326203& 4658099\\
\, & \, & 79794157&   1607252663& 37447368857& 995564440951\\
\, & \, & 29858589333061...\\
    86&  443& 3109& 28573& 332099& 4748047\\
 \, & \, &  81428323&  1641908027& 38291437141& 1018893116299\\
 \, & \, & 30582699050611...\\
   87&  449& 3169& 29153& 339601& 4863959\\
   \, & \, & 83543071& 1686826109& 39386748617& 1049194449883\\
   \, & \, & 31524064728311...\\
 88&  457& 3229& 29803& 347849& 4989697\\
 \, & \, & 85839547&  1735649329& 40578571003& 1082201297941\\
 \, & \, & 32550506359429...\\
   90&  463& 3299& 30557& 357473& 5138719\\
   \, & \, & 88565483&  1793681753& 41997140089& 1121535591721\\
   \, & \, & 33775078562347...\\
    91&  467& 3319& 30781& 360293& 5182717\\
    \, & \, & 89369047& 1810798861& 42415879469& 1133155938589\\
    \, & \, & 34137123380603...\\
    92&  479& 3407& 31667& 371981& 5363167\\
    \, & \, & 92678347&  1881428537& 44145738083& 1181205761389\\
    \, & \, & 35635464099689...\\
\end{array}
$$

{\bf Appendix 1}\\

{\bf Matrix $\mathbf{[\ov{M}\phantom{d} ^2P]}$}\ldots\quad Continuation 5.
$$
\begin{array}{llllll}
\hline\\
  93&  487&  3469& 32341& 380557& 5496349\\
   \, & \, & 95121911& 1933651711& 45426482839& 1216826411041\\
   \, & \, & 36747532444747...\\
    94&  491& 3517& 32797& 386401& 5587537\\
    \, & \, & 96797411&  1969496239& 46306458839& 1241322670799\\
    \, & \, & 37512927359291...\\
   95&  499&  3559& 33203& 391711& 5670851\\
   \, & \, & 98330021&  2002298621& 47112340151& 1263771327193\\
   \, & \, & 38214783465337...\\
    96&  503& 3593& 33569& 396269& 5741453\\
    \, & \, & 99630571&  2030158657& 47797243919& 1282861540019\\
    \, & \, & 38811965770483...\\
    98&  521& 3733& 35023& 415253& 6037513\\
    \, & \, & 105089261&  2147305243& 50681376121& 1363360331743\\
    \, & \, & 41333311232987...\\
    99&  523& 3761& 35311& 418961& 6095731\\
    \, & \, & 106166089&  2170447637& 51251887327& 1379303865481\\
    \, & \, & 41833278300773...\\
   100&  541& 3911& 36887& 439357& 6415081\\
   \, & \, & 112073683&  2297602183& 54391267121& 1467155677657\\
   \, & \, & 44591559921641...\\
   102&  557&  4027&  38153&  455849&  6673993\\
   \, & \, &  116881321&  2401362767& 56958606937&  1539140110927\\
   \, & \, &  46855727983837...\\
   104&  569&  4133&  39239&  470207&  6898807\\
   \, & \, &  121064467&  2491797367& 59200082443&  1602086508713\\
   \, & \, &  48838469899327...\\
 105&  571&  4153&  39451&  472837&  6940103\\
 \, & \, &  121834483&  2508461203&  59613478459&  1613705610163\\
 \, & \, & 49204743622123...\\
   106&  577&   4217&  40151&  481847&  7081709\\
   \, & \, &  124469621& 2565499711& 61029312569&  1653521623993\\
   \, & \, & 50460527025823...\\
   108&  593&  4339&  41491&  499403&  7359427\\
   \, & \, &  129647857&  2677808011&  63821022049&  1732128413677\\
   \, & \, &  52942646093899...\\
   110&  601&  4421&  42293&  510031&  7528669\\
   \, & \, &  132814411&  2746597487&  65533394977&  1780407360517\\
   \, & \, &  54468962620717...\\
   111&  607&  4463&  42697&  515401&  7612799\\
   \, & \, &  134389627&  2780844971&  66386576369&  1804479121591\\
   \, & \, &  55230488801623...\\
\end{array}
$$

{\bf Appendix 1}\\

{\bf Matrix $\mathbf{[\ov{M}\phantom{d} ^2P]}$}\ldots\quad Continuation 6.
$$
\begin{array}{llllll}
\hline\\
 112&  613&  4517&  43283&  522829&  7730539\\
   \, & \, &  136593931&  2828789699&  67581794939&  1838220650251\\
   \, & \, &  56298481067219...\\

114& 619& 4567& 43889& 530773&  7856939\ldots\\
115& 631& 4663& 44879& 543967&  8066533\ldots\\
116& 641& 4759& 45971& 558643& 8300687\ldots\\
117& 643& 4787& 46279& 562711& 8365481\ldots\\
118& 647& 4801& 46451& 565069& 8402833\ldots\\
119& 653& 4877& 47297& 576203&  8580151\ldots\\
120& 659& 4933& 47857& 583523& 8696917\ldots\\
121& 661& 4943& 47963& 584999& 8720227\ldots\\
122& 673& 5021& 48821& 596243& 8900383\ldots\\
123& 677& 5059& 49207& 601397& 8982923\ldots\\
124& 683& 5107& 49739& 608459& 9096533\ldots\\
125& 691& 5189& 50591& 619739& 9276991\ldots\\
126& 701& 5281& 51599& 633467& 9498161\ldots\\
127& 709& 5381& 52711& 648391& 9737333\ldots\\
128& 719& 5441& 53353& 657121& 9878657\ldots\\
  129& 727& 5503& 54013& 665843& 10020343\ldots\\
  130& 733& 5557& 54601& 673793&   10147877\ldots\\
    132& 743& 5651& 55681& 688249& 10382033\ldots\\
  133& 751& 5701& 56197& 695239&    10493953\ldots\\
  134& 757& 5749& 56701& 702173& 10606223\ldots\\
  135& 761& 5801& 57193& 708479& 10707449\ldots\\
   136& 769& 5851& 57751& 715969&  10829519\ldots\\
     138& 787& 6037& 59723& 742681& 11261903\ldots\\
      140& 809& 6217& 61819& 771079& 11723507\ldots\\
  141& 811& 6229& 61979& 773317& 11760029\ldots\\
  142& 821& 6311& 62921& 786053&    11967047\ldots\\
  143& 823& 6323& 63059& 788009& 11999111\ldots\\
  144& 827& 6353& 63391& 792413& 12071197\ldots\\
  145& 829& 6361& 63467& 793511&    12089177\ldots\\
  146& 839& 6469& 64679& 809917& 12356863\ldots\\
   147& 853& 6599& 66089& 828923& 12667463\ldots\\
   148& 857& 6653& 66749& 838091&    12816389\ldots\\
       150& 863& 6691& 67157& 843613& 12907091\ldots\\
        152& 881& 6841& 68821& 866329& 13280819\ldots\\
  153& 883& 6863& 69109& 870161& 13343881\ldots\\
  154& 887& 6899& 69491& 875519&    13431967\ldots\\
  155& 907& 7057& 71287& 900157& 13836751\ldots\\
   \cdot& \cdot& \cdot& \cdot& \cdot& \cdot
\end{array}
$$

}

\end{document}